\documentclass[reqno]{amsart}

\usepackage{amssymb}
\usepackage[T1]{fontenc}
\usepackage{graphicx}
\usepackage{ifpdf}
\usepackage{xypic}

\author{Konrad Sch\"obel}
\title
	[$\PU(2)$-Instantons on Class $\VII_0$ Surfaces with $b_2=1$]
	{Moduli Spaces of $\PU(2)$-Instantons \\ on Minimal Class~VII Surfaces with $b_2=1$}
\email{schoebel@cmi.univ-mrs.fr}
\address{
	Laboratoire d'Analyse, Topologie et Probabilit\'es \\
	Centre de Math\'ematiques et Informatique \\
	Universit\'e de Provence \\
	39, rue F. Joliot Curie \\
	13\,453 Marseille Cedex~13 \\
	France
}
\thanks{E-mail: \texttt{schoebel@cmi.univ-mrs.fr}}

\ifpdf\fi

\numberwithin{equation}{section}

\newtheorem{theorem}{Theorem}[section]
\newtheorem{proposition}[theorem]{Proposition}
\newtheorem{lemma}[theorem]{Lemma}
\newtheorem{corollary}[theorem]{Corollary}
\newtheorem{definition}[theorem]{Definition}
\newtheorem{remark}[theorem]{Remark}

\newtheorem{examples}[theorem]{Examples}
\newtheorem{notation}[theorem]{Notation}

\DeclareSymbolFont{rsfs}{U}{rsfs}{m}{n}
\DeclareSymbolFontAlphabet{\mathscr}{rsfs}

\newcommand{\name}[1]{\textsc{#1}}
\newcommand{\dfn}[1]{\emph{#1}}

\renewcommand{\deg}{\operatorname{deg}}
\renewcommand{\emptyset}{\varnothing}
\renewcommand{\ge}{\geqslant}
\renewcommand{\le}{\leqslant}
\newcommand{\lle}{(\le)}

\newcommand{\comma}{\,,}
\newcommand{\fullstop}{\,.}
\newcommand{\coloneq}{:=}
\newcommand{\eqcolon}{=:}
\newcommand{\isom}{\cong}

\newcommand{\Z}{\mathbb{Z}}
\newcommand{\N}{\mathbb{N}}
\newcommand{\R}{\mathbb{R}}
\newcommand{\C}{\mathbb{C}}

\newcommand{\sheaf}[1]{\mathscr{#1}}
\newcommand{\hol}[1]{\mathcal{#1}}

\newcommand{\VII}{{\rm VII} }
\newcommand{\dual}{\vee}

\DeclareMathOperator{\id}{id}
\DeclareMathOperator{\vol}{vol}
\DeclareMathOperator{\Ext}{Ext}

\DeclareMathOperator{\Pic}{Pic}
\DeclareMathOperator{\GL}{GL}
\DeclareMathOperator{\U}{U}
\DeclareMathOperator{\SU}{SU}
\DeclareMathOperator{\PU}{PU}
\DeclareMathOperator{\Ad}{Ad}
\DeclareMathOperator{\ad}{ad}

\newcommand{\roots}{R}
\newcommand{\queer}{Q}
\newcommand{\unstable}{U}

\begin{document}

\begin{abstract}
	
	We describe explicitly the moduli spaces $\mathcal M^{\rm pst}_g(S,E)$ of polystable holomorphic structures $\hol E$ with
	$\det\hol E\isom\hol K$ on a rank $2$ vector bundle $E$ with $c_1(E)=c_1(K)$ and $c_2(E)=0$ for all minimal class \VII
	surfaces $S$ with $b_2(S)=1$ and with respect to all possible \name{Gauduchon} metrics $g$.  These surfaces $S$ are
	non-elliptic and non-\name{K\"ahler} complex surfaces and have recently been completely classified \cite{Teleman}.  When $S$
	is a half or parabolic \name{Inoue} surface, $\mathcal M^{\rm pst}_g(S,E)$ is always a compact one-dimensional complex
	disc.  When $S$ is an \name{Enoki} surface, one obtains a complex disc with finitely many transverse self-intersections
	whose number becomes arbitrarily large when $g$ varies in the space of \name{Gauduchon} metrics.  $\mathcal
	M^{\rm pst}_g(S,E)$ can be identified with a moduli space of $\PU(2)$-instantons.  The moduli spaces of simple bundles of
	the above type leads to interesting examples of non-\name{Hausdorff} singular one-dimensional complex spaces.
	
\end{abstract}

\maketitle

\bigskip

\section{Introduction}

In gauge theory, moduli spaces of anti-self-dual connections have led to striking results in differential four-manifold
geometry; they are the main tools in the construction of the \name{Donaldson} polynomial invariants.  However, the explicit
computation of these moduli spaces in concrete situations is in general very difficult.  On the other hand, when the base
manifold is a complex surface, the \name{Kobayashi}-\name{Hitchin} correspondence establishes a real analytic isomorphism
between the moduli spaces of (irreducible) anti-self-dual connections and polystable (stable) holomorphic structures on a fixed
differentiable vector bundle and makes thus possible the application of complex geometric methods for the computation of
gauge-theoretical moduli spaces.  \name{S.~K.~Donaldson} gave the first complete proof of this relationship on algebraic
surfaces and used it to explicitly compute moduli spaces and the corresponding invariants for \name{Dolgachev} surfaces.  This
led to the first example of pairs of homeomorphic but not diffeomorphic four-manifolds \cite{Donaldson}.

Subsequently this strategy was carried out for a large variety of algebraic surfaces \cite{Okonek&Van_de_Ven86, Buchdahl87,
Friedman, Kotschick, Okonek&Van_de_Ven89, Friedman&Morgan, Donagi, Friedman&Morgan&Witten}.  However, it becomes very hard for
non-algebraic surfaces due to the presence of non-filtrable holomorphic bundles in the moduli space.  A complete classification
of such bundles is considered to be an extremely difficult problem on non-elliptic surfaces because of the lack of a general
method of construction and parametrisation.  On the other hand, on elliptic surfaces it could be solved for a number of
non-\name{K\"ahler}ian elliptic surfaces \cite{Braam&Hurtubise, Luebke&Teleman, Teleman98, Toma01, Moraru, Brinzanescu&Moraru}.
Note that for elliptic fibrations one solves this problem by regarding the restrictions to the fibres which (generically) are
elliptic curves on which the classification of holomorphic bundles is well understood \cite{Atiyah&Bott, Potier&Verdier}.  This
strategy is called the graph method and was used by \name{P.~J.~Braam} and \name{J.~Hurtubise} to obtain the first explicit
example of an $\SU(2)$-instanton moduli space on a non-\name{K\"ahler} surface, namely an elliptic \name{Hopf} surface
\cite{Braam&Hurtubise}.  In this article we now compute moduli spaces of holomorphic bundles on \emph{all} minimal class \VII
surfaces with $b_2=1$, endowed with \emph{all} possible \name{Gauduchon} metrics.  Being the first example of moduli spaces on
surfaces that are both non-\name{K\"ahler} and non-elliptic, this is the reason why one expects essentially new phenomena for
the behaviour of moduli spaces in general.

Our method to overcome the main difficulty of controlling non-filtrable bundles is the following:  We first classify filtrable
bundles and then show, using gauge-theory, that only a particular non-filtrable bundle can exist.  In \cite{Teleman} it was
shown that the moduli space does not contain a compact component consisting of both filtrable and non-filtrable bundles.  We
then show that the moduli space does not contain any compact component at all.  This is true on blown-up primary \name{Hopf}
surfaces by a recent result of \name{M.~Toma} \cite{Toma06} and we conclude using a deformation argument, since any minimal
class \VII surface surface containing a global spherical shell (see below) is the degeneration of a blown-up primary \name{Hopf}
surface \cite{Kato78}.

A first interesting property of our moduli spaces is that the filtrable bundles are generic.  This is surprising, because on
\name{K\"ahler} surfaces the filtrable locus is a countable union of \name{Zariski}-closed sets and also in all formerly known
examples on non-\name{K\"ahler} surfaces it was found to be \name{Zariski}-closed.

Class \VII surfaces with $b_2=1$ are of particular interest in the light of the classification problem of complex surfaces.  In
the early 1960ies, \name{K.~Kodaira} classified connected compact complex surfaces (\dfn{surfaces} for short) into seven
classes \cite{Kodaira-bI}.  Six of them are quite well understood but the seventh \cite{Kodaira-bII} has resisted a complete
classification until the present day.  A surface $S$ is said to be of class \VII if it has \name{Kodaira} dimension
$\operatorname{kod}(S)=-\infty$ and first \name{Betti} number $b_1(S)=1$.  It can be blown down to a unique minimal model,
i.\,e.\ a unique class \VII surface not being the blow-up of another one.  We denote the subclass of minimal class \VII surfaces
by $\VII_0$.  Class $\VII_0$ surfaces with second \name{Betti} number $b_2=0$ are classified:  They are either \name{Hopf} or
\name{Inoue} surfaces \cite{Bogomolov82, Li&Yau&Zheng94, Teleman94}.  As to class $\VII_0$ surfaces with $b_2>0$, all known
examples admit a so-called global spherical shell and can be explicitly constructed by successive blow-ups and holomorphic
surgery \cite{Kato78}.  On the other hand, every class $\VII_0$ surface $S$ with exactly $b_2(S)$ rational curves possesses a
global spherical shell \cite{Dloussky&Oeljeklaus&Toma}.  The global spherical shell conjecture now states that every class
$\VII_0$ surface has such a global spherical shell and would reduce the classification of class \VII surfaces to finding
sufficiently many curves.

This was recently done by \name{A.~Teleman} for the subclass $\VII_0^1$ of class $\VII_0$ surfaces with $b_2=1$ \cite{Teleman}.
Supposing there did not exist any complex curves on the surface he constructed a contradiction for the moduli space of
polystable holomorphic structures $\hol E$ with $\det\hol E\isom\hol K$ on a fixed complex vector bundle $E$ with
$c_1(E)=c_1(K)$ and $c_2(E)=0$.  By the above, this accomplishes the classification of class $\VII_0^1$ surfaces:  Each class
$\VII_0^1$ surface is biholomorphic to either the half \name{Inoue} surface \cite{Inoue77}, the parabolic \name{Inoue} surface
\cite{Inoue74} or an \name{Enoki} surface \cite{Enoki80}.  We in turn now compute explicitly this moduli space for each of these
surfaces and describe its properties in detail.  This is possible with respect to any \name{Gauduchon} metric, due to a
recent result classifying the possible degree maps on non-\name{K\"ahler} surfaces \cite{Buchdahl00,Teleman05}.  We finally
remark that the methods used can be extended to show the existence of a curve in the case $b_2=2$ \cite{Teleman06b}.

The expected complex dimension of the above moduli space is
\begin{equation}
	\label{eq:dim}
	-\chi(\hol End_0\hol E)=\big(4c_2(E)-c_1(E)^2\big)-\tfrac32\big(b_2^+(S)-b_1(S)+1\big)=1
	\comma
\end{equation}
but there are two deeper reasons for this particular choice of the \name{Chern} classes of $E$.  Firstly, it allows one to write
filtrable holomorphic bundles $\hol E$ as extensions of certain holomorphic line bundles.  Secondly, it assures that the moduli
space of anti-self-dual connections on $E$ is compact so that the moduli space of stable holomorphic structures on $E$, embedded
via the \name{Kobayashi}-\name{Hitchin} correspondence, can be compactified by adding only the irreducible part.  This
compactification is crucial in the step determining possible non-filtrable bundles.

The moduli spaces we get are compact one-dimensional complex discs when the surface is a half or parabolic \name{Inoue} surface.
In the ``generic'' case of an \name{Enoki} surface it is a compact one-dimensional complex disc too, but with finitely many
transverse self-intersections.  The number of these singularities is unbounded when the metric varies in the space of
\name{Gauduchon} metrics.  This shows that there are infinitely many homeomorphism types of moduli spaces although there are
only finitely many topological splittings of the underlying vector bundle.  Furthermore, having a boundary, these moduli spaces
are not complex spaces.  This is in contrast to algebraic surfaces, where the \name{Uhlenbeck} compactification is known to be 
an algebraic variety \cite{Li}, and to all known examples on non-algebraic surfaces.  It will be one of our next steps to study
the behaviour of the natural \name{Hermit}ian metric \cite{Luebke&Teleman} near this boundary.

Let us finally point out that our results could only be obtained through a close interplay between complex geometry and gauge
theory.  Although nowadays \name{Seiberg}-\name{Witten} theory has widely replaced \name{Donaldson} theory, recent developments
show that \name{Donaldson} theory on definite $4$-manifolds with $b_1\ge1$ is still an interesting open subject
\cite{Teleman06a}.

\medskip

The structure of this article is the following:  In the next section we briefly review the necessary properties of class
$\VII_0^1$ surfaces and summarise their classification.  Then we parametrise filtrable holomorphic bundles in the moduli spaces
(section~\ref{sec:filtrable}), examine its local structure (section~\ref{sec:local}) and the stability condition
(section~\ref{sec:stable}).  In section~\ref{sec:boundary} we give the boundary structure of the moduli spaces of polystable
bundles.  Finally we determine non-filtrable bundles (section~\ref{sec:centres}) which leads to a complete description of the
entire moduli spaces in the last section.

\medskip

\textbf{Acknowledgements}:
I would like to thank my PhD supervisor \name{Andrei Teleman} for his constant help and numerous hints and comments.
Furthermore I would like to thank \name{Matei Toma} for providing his results as well as \name{Karl Oeljeklaus} and \name{George
Dloussky} for useful discussions on the subject.

\section{Minimal class \VII surfaces with $b_2=1$}
\label{sec:VII-surfaces}

Let $S$ be a class \VII surface, i.\,e.\ a compact complex surface with first \name{Betti} number $b_1(S)=1$ and \name{Kodaira}
dimension $\operatorname{kod}(S)=-\infty$.  By definition, the condition on the \name{Kodaira} dimension means that tensor
powers of the canonical holomorphic line bundle $\hol K$ do not admit any non-trivial holomorphic sections, i.\,e.\ $H^0(\hol
K^{\otimes n})=0$ for $n\ge1$.  For such a surface the \name{Chern} classes are given by \cite{Kodaira-bI}
\begin{equation}
	\label{eq:Chern_classes}
	c_2(S)=-c_1(S)^2=b_2(S)
	\fullstop
\end{equation}

Suppose now that $S$ is of class $\VII_0^1$, i.\,e.\ minimal with second \name{Betti} number $b_2(S)=1$.  As mentioned in the
introduction, \name{Teleman} proved that in this case there exists at least one complex curve on $S$ \cite{Teleman}.  But any
class $\VII_0^1$ surface containing a curve is biholomorphic to one of the following surfaces \cite{Nakamura-I}:
\begin{itemize}
	\item
		A \dfn{half \name{Inoue} surface} \cite{Inoue77}.  It contains only a single complex curve, namely a singular rational
		curve $C$ with one node and self intersection $-1$.  The canonical bundle is given by
		\begin{equation}
			\label{eq:half}
			\hol K=\hol F\otimes\hol O(-C)
		\end{equation}
		where $\hol F$ is the unique non-trivial square-root of the trivial holomorphic line bundle $\hol O$ (see below).
		We have
		\[
			c_1(\hol O(C))=-c_1(\hol K)
			\fullstop
		\]
	\item
		A surface in the family studied by \name{Enoki} \cite{Enoki80} containing only a single complex curve, namely a singular
		rational curve $C$ with one node and self intersection $0$.  We have
		\[
			c_1(\hol O(C))=0
			\fullstop
		\]
		There is no expression for the canonical bundle $\hol K$ as in the other two cases.  We will refer to class $\VII_0^1$
		surfaces of this type as \dfn{\name{Enoki} surfaces}.
	\item A \dfn{parabolic \name{Inoue} surface} \cite{Inoue74}.  It contains precisely two complex curves, namely a singular
		rational curve $C$ with one node and self intersection $0$ and an elliptic curve $E$ with self-intersection
		$-1$.  Both curves are disjoint.  The canonical bundle is given by
		\begin{equation}
			\label{eq:parabolic}
			\hol K=\hol O(-C-E)
			\fullstop
		\end{equation}
		We have
		\[
			c_1(\hol O(C))=0
			\qquad
			c_1(\hol O(E))=-c_1(\hol K)
			\fullstop
		\]
\end{itemize}
The \name{Chern} classes above follow from the intersection numbers since $H^2(S,\Z)$ is torsion free (see below).

\begin{notation}
	The family of class $\VII_0^1$ surfaces constructed and classified by \name{Enoki} \cite{Enoki80, Enoki81} is characterised
	by the existence of a divisor $D>0$ with $D^2=0$.  As such it includes the parabolic \name{Inoue} surface.  Nevertheless, to
	simplify the exposition we agree that \textbf{in this article we do not consider the parabolic \name{Inoue} surface as an
	\name{Enoki} surface}.
	
	Unless otherwise stated, $S$ will always denote a class $\VII_0^1$ surface, i.\,e.\ one of the three types above.
\end{notation}

\begin{remark}
	As a two parameter family \name{Enoki} surfaces represent the generic case of class $\VII_0^1$ surfaces.  The half and the
	parabolic \name{Inoue} surface appear as degenerations of them.
\end{remark}

The existence of a rational curve on a class $\VII_0^1$ surface implies the existence of a so-called global spherical shell
\cite{Dloussky&Oeljeklaus&Toma}.  Surfaces admitting a global spherical shell can be constructed by successive blow-ups of the
unit ball in $\C^2$ and a subsequent holomorphic surgery \cite{Kato78, Dloussky84}.  A consequence of this construction is that
all such surfaces are degenerations of blown-up primary \name{Hopf} surfaces.  In particular they are all diffeomorphic with
fundamental group $\pi_1(S)\isom\Z$.  Thus $H_1(S,\Z)\isom\Z$ is free and from the universal coefficient theorem we conclude
$H^2(S,\Z)\isom\Z$ because $b_2(S)=1$.  Furthermore, from \eqref{eq:Chern_classes} we see that $c_1(K)^2=-1$, showing that
$c_1(K)$ is a generator of $H^2(S,\Z)$.

In the following we will frequently use the correspondence between line bundle morphisms $\hol M_1\to\hol M_2$ and the sections
of $\hol M_1^\dual\otimes\hol M_2$ they define.  In particular every such morphism is the zero morphism if the corresponding
bundle does not admit non-trivial sections, i.\,e.\ if $H^0(\hol M_1^\dual\otimes\hol M_2)=0$.
\begin{remark}
	\label{rem:sections}
	Note that a line bundle admits non-trivial sections if and only if it is isomorphic to $\hol O(D)$ for a divisor $D\ge0$ on
	$S$, i.\,e.\ if it is of the form $\hol O(rC)$ on the half \name{Inoue} or an \name{Enoki} surface and $\hol O(rC+sE)$ on
	the parabolic \name{Inoue} surface for some $r,s\in\N$.  This shows in particular that line bundles $\hol M$ on class
	$VII_0^1$ surfaces with $c_1(\hol M)=c_1(\hol K^{\otimes n})$ do not admit non-trivial sections if $n\ge1$, a fact we will
	use frequently below without further mention.

	The divisor $D$ is the zero divisor of a section in the line bundle and uniquely determined since class $\VII_0^1$ surfaces
	do not admit non-constant meromorphic functions.  In particular we have $\dim H^0(\hol M)\le1$ for line bundles $\hol M$ on
	class $\VII_0^1$ surfaces and if $\hol M$ is non-trivial then either $\hol M$ or $\hol M^\dual$ does not admit non-trivial
	sections.
\end{remark}

The exponential sequence $0\to\Z\to\hol O\xrightarrow{\exp}\hol O^*\to0$ gives rise to the long exact cohomology sequence
\[
	\ldots\to H^1(S,\Z)\to H^1(S,\hol O)\xrightarrow{\exp^1}H^1(S,\hol O^*)\xrightarrow{c_1}H^2(S,\Z)\to\ldots
	\fullstop
\]
Here $\Pic(S)\coloneq H^1(S,\hol O^*)$ is the \name{Picard} group, the \name{Abel}ian group of isomorphism classes of
holomorphic line bundles on $S$ with group multiplication induced by the tensor product.  On the other hand $H^2(S,\Z)$
classifies isomorphism classes of complex line bundles via the first \name{Chern} class.  The connecting operator is just the
group homomorphism that associates to a holomorphic line bundle the first \name{Chern} class of its underlying topological line
bundle.  Its kernel, the image of $\exp^1$, is the subgroup $\Pic^0(S)$ of holomorphic structures on the topologically trivial
line bundle.  The \name{Picard} group $\Pic(S)$ has the structure of a complex \name{Lie} group and $\exp^1$ is an \'etale
morphism \cite{Luebke&Teleman}.

Since $H^1(S,\Z)$ is torsion free and $b_1(S)=1$ we have $H^1(S,\Z)\isom\Z$.  On the other hand, on class \VII surfaces the
natural inclusion $\C\hookrightarrow\hol O$ induces an isomorphism $H^1(S,\C)\stackrel{\isom}{\longrightarrow}H^1(S,\hol O)$
\cite{Kodaira-bI} showing that for $b_1(S)=1$ there is a group isomorphism
\[
	\Pic^0(S)\isom\C^*
	\fullstop
\]
In particular, every holomorphic line bundle in $\Pic^0(S)$ has exactly two roots in $\Pic^0(S)$ which differ by the non-trivial
root of $\hol O$ which we will denote by $\hol F$:
\[
	\hol F\otimes\hol F=\hol O
	\qquad
	\hol F\not\isom\hol O
	\fullstop
\]
Remark that in contrast to \name{K\"ahler} surfaces $\Pic^0(S)$ is non-compact here.

\section{Filtrable holomorphic bundles}
\label{sec:filtrable}

On surfaces, topological complex vector bundles are classified up to isomorphisms by their rank and their first two \name{Chern}
classes.  We fix once and for all a complex vector bundle $E$ on $S$ with
\begin{subequations}
	\label{eq:fix}
	\begin{equation}
		\label{eq:fix_E}
		\operatorname{rank}E=2
		\qquad
		c_1(E)=c_1(K)
		\qquad
		c_2(E)=0
		\comma
	\end{equation}
where $K$ is the canonical complex line bundle.  Since $c_1(\det E)=c_1(E)$ this implies $\det E\isom K$.  In the following we
will study the simple holomorphic structures $\hol E$ on $E$ with determinant
	\begin{equation}
		\label{eq:fix_det}
		\det\hol E\isom\hol K
		\fullstop
	\end{equation}
\end{subequations}

At first we investigate \emph{filtrable} bundles of type \eqref{eq:fix}, because they admit a relatively simple description as
extensions of certain holomorphic line bundles.  In general a rank two bundle is filtrable if it admits a rank one subsheaf,
but the notion simplifies considerably for surfaces:
\begin{definition}
	\label{def:filtrable}
	A holomorphic rank two vector bundle $\hol E$ on a complex surface $S$ is \dfn{filtrable} if and only if one of the
	following equivalent conditions is satisfied:
	\begin{enumerate}
		\item \label{def:filtrable:1}
			$\hol E$ has a rank one subsheaf $\sheaf S$.
		\item \label{def:filtrable:2}
			$\hol E$ has a locally free rank one subsheaf $\hol L$.
		\item \label{def:filtrable:3}
			There exist holomorphic line bundles $\hol L$ and $\hol R$ on $S$ that fit into a short exact sequence of the form
			\begin{equation}
				\label{eq:filtrable:3}
				0\quad\longrightarrow\quad\hol L\quad\longrightarrow\quad\hol E\quad\longrightarrow\quad\hol R\otimes\sheaf I_Z\quad\longrightarrow\quad0
			\end{equation}
			where $\sheaf I_Z$ is the ideal sheaf of a dimension zero locally complete intersection $Z\subset S$.
	\end{enumerate}
\end{definition}
The proof of the equivalence is standard, see for example \cite{Elencwajg&Forster}.

The first reason for choosing $E$ to satisfy \eqref{eq:fix_E} is that in this case we get rid of the (possibly very complicated)
ideal sheaf $\sheaf I_Z$ in \eqref{eq:filtrable:3}:

\begin{proposition}
	\label{prop:ZLR}
	On a class $\VII_0^1$ surface $S$ we have $Z=\emptyset$ and either $c_1(\hol L)=0$ or $c_1(\hol R)=0$ in
	\eqref{eq:filtrable:3} under the assumption \eqref{eq:fix_E}.
\end{proposition}

\begin{proof}
	Since $c_1(\hol K)$ is a generator of $H^2(S,\Z)\isom\Z$ we set $c_1(\hol L)=n\cdot c_1(\hol K)$ with $n\in\Z$.  Computation
	of the \name{Chern} classes of $\hol E=(\hol E\otimes\hol L^\dual)\otimes\hol L$ yields, since $c_1(\hol K)^2=-1$, 
	\[
		\begin{array}{r@{\;}c@{\;}c@{\;}c@{\;}l}
			c_1(\hol K)
			&=&c_1(\hol E)
			&=&c_1(\hol E\otimes\hol L^\dual)+2c_1(\hol L)
				\qquad\text{and}\\
			0
			&=&c_2(\hol E)
			&=&c_2(\hol E\otimes\hol L^\dual)+c_1(\hol E\otimes\hol L^\dual)c_1(\hol L)+c_1(\hol L)^2
			 =|Z|+n(n-1)
			\fullstop
		\end{array}
	\]
	Here $|Z|$ denotes the number of points in $Z$, counted with multiplicities.  But the last equality can only be satisfied
	if $|Z|=0$, i.\,e.\ $Z=\emptyset$, and $n=0$ or $1$.
\end{proof}

Now note that the determinant of a the central term in a line bundle extension is the tensor product of the two corresponding
line bundles.

\begin{corollary}
	\label{cor:extensions}
	Any filtrable holomorphic vector bundle $\hol E$ of type \eqref{eq:fix} on a class $\VII_0^1$ surface is the central term of
	an extension of one of the following two types
	\begin{subequations}
		\label{eq:EA}
		\begin{alignat}{8}
			\label{eq:E}
			0&\quad\longrightarrow&&\hol L&&&\longrightarrow\quad&\hol E&\quad&\to\quad&\hol L^\dual&\;\otimes\;&\hol K&&\quad\to\quad&0
		\\[\medskipamount]
			\label{eq:A}
			0&\quad\to\quad&\hol R^\dual&\;\otimes\;&\hol K&&\quad\to\quad&\hol E&\quad&\longrightarrow&&\hol R&&&\longrightarrow\quad&0
		\end{alignat}
	\end{subequations}
	where $\hol L,\hol R\in\Pic^0(S)$.
	
	Moreover, given a line bundle inclusion $\hol M\hookrightarrow\hol E$ into a bundle $\hol E$ of type \eqref{eq:fix}, we have
	either $c_1(\hol M)=0$ and the inclusion extends to \eqref{eq:E} with $\hol L\isom\hol M$ or we have $c_1(\hol M)=c_1(\hol
	K)$ and it extends to \eqref{eq:A} with $\hol R^\dual\otimes\hol K\isom\hol M$.
\end{corollary}

The following lemma shows that the existence of non-trivial extensions \eqref{eq:EA} and the uniqueness of their central terms
is determined by the existence of sections in certain line bundles.  Line bundle extension $0\to\hol M_1\to\hol E\to\hol
M_2\to0$ or equivalently $0\to\hol M_1\otimes\hol M_2^\dual\to\hol E\otimes\hol M_2^\dual\to\hol O\to0$ are determined by the
image of the constant $1$ section in $\hol O$ under the connecting operator $H^0(\hol O)\to H^1(\hol M_1\otimes\hol M_2^\dual)$
in the associated cohomology sequence and vice versa.  In particular extensions which differ by a non-zero constant in the
classifying space $\Ext^1(\hol M_2,\hol M_1)\coloneq H^1(\hol M_1\otimes\hol M_2^\dual)$ have isomorphic central terms.  To
compute the dimension of these spaces we will use the \name{Riemann}-\name{Roch} theorem which, using \eqref{eq:Chern_classes}
and combined with the \name{Serre} duality, takes the particular form
\begin{equation}
	\label{eq:RRSD}
	h^0(\hol M)-h^1(\hol M)+h^0(\hol M^\dual\otimes\hol K)
	=
	\tfrac12\:c_1(\hol M)\:\big(c_1(\hol M)-c_1(\hol K)\big)
\end{equation}
for a holomorphic line bundle $\hol M$ on $S$, where $h^p(\hol M)\coloneq\dim H^p(S,\hol M)$ \cite{BHPvV}.  To simplify the
notation we write $\hol L^2$ and $\hol L^{-2}$ for $\hol L\otimes\hol L$ and $\hol L^\vee\otimes\hol L^\vee$ respectively.

\begin{proposition}
	\label{prop:param9n}
	\begin{enumerate}
		\item \label{prop:param9n:1}
			For every holomorphic line bundle $\hol L\in\Pic^0(S)\setminus\queer(S)$, where
			\[
				\queer(S)\coloneq\{\hol L\in\Pic^0(S)\colon H^0(\hol L^2\otimes\hol K^\dual)\not=0\}
				\comma
			\]
			there is a non-trivial extension
			\begin{subequations}
			\label{eq:E_L&A_R}
			\begin{equation}
				\label{eq:E_L}
				0\quad\longrightarrow\quad\hol L\quad\longrightarrow\quad\hol E_{\hol L}\quad\longrightarrow\quad\hol L^\dual\otimes\hol K\quad\longrightarrow\quad0
			\end{equation}
			with an (up to isomorphisms) uniquely determined central term $\hol E_{\hol L}$.
			If $\hol L\in\queer(S)$ then the isomorphism classes of central terms in non-trivial extensions of the form \eqref{eq:E}
			are parametrised by $\C\rm P^1$.
		\item \label{prop:param9n:2}
			For every holomorphic line bundle $\hol R\in\roots(S)$,
			where
			\[
				\roots(S)\coloneq\{\hol R\in\Pic^0(S)\colon H^0(\hol R^2)\not=0\}
				\comma
			\]
			there is a non-trivial extension
			\begin{equation}
				\label{eq:A_R}
				0\quad\longrightarrow\quad\hol R^\dual\otimes\hol K\quad\longrightarrow\quad\hol A_{\hol R}\quad\longrightarrow\quad\hol R\quad\longrightarrow\quad0
			\end{equation}
			\end{subequations}
			with an (up to isomorphisms) uniquely determined central term $\hol A_{\hol R}$.
			If $\hol R\in\Pic^0(S)\setminus\roots(S)$ there are no non-trivial extensions of the form \eqref{eq:A}.
	\end{enumerate}
\end{proposition}

\begin{proof}
	Extensions of type \eqref{eq:A_R} are classified by $\Ext^1(\hol R,\hol R^\dual\otimes\hol K)\isom H^1(\hol
	R^{-2}\otimes\hol K)$.  From formula \eqref{eq:RRSD} for $\hol M=\hol R^{-2}\otimes\hol K$ we obtain $\dim\Ext^1(\hol R,\hol
	R^\dual\otimes\hol K)=h^0(\hol R^2)$ since $H^0(\hol R^{-2}\otimes\hol K)=0$.  This proves (\ref{prop:param9n:2}) because
	$h^0(\hol R^2)=0$ or $1$.
	Likewise, extensions of type \eqref{eq:E_L} are classified by $\Ext^1(\hol L^\dual\otimes\hol K,\hol L)\isom H^1(\hol
	L^2\otimes\hol K^\dual)$ and from formula \eqref{eq:RRSD} we obtain $\dim\Ext^1(\hol L^\dual\otimes\hol K,\hol L)=1+h^0(\hol
	L^2\otimes\hol K^\dual)$ since $H^0(\hol L^{-2}\otimes\hol K^2)=0$.  This proves the first part of (\ref{prop:param9n:1}).

	In the case $\hol L\in\queer(S)$ we have $\dim\Ext^1(\hol L^\dual\otimes\hol K,\hol L)=1+h^0(\hol L^2\otimes\hol K^\dual)=2$
	because $0\not=h^0(\hol L^2\otimes\hol K^\dual)\le1$.  Let $\varphi\colon\hol E_1\to\hol E_2$ be a bundle isomorphism between
	the central terms of two different extensions in the following diagram:
	\begin{equation}
		\label{diag:E12}
		\xymatrix{
			0\ar[r]&\hol L\ar^{\alpha_1}[r]\ar@{-->}[d]&\hol E_1\ar[r]\ar_\varphi[d] &\hol L^\dual\otimes\hol K\ar[r]\ar@<4pt>@{-->}[d]&0\\
			0\ar[r]&\hol L\ar[r]                       &\hol E_2\ar^{\beta_2\quad}[r]&\hol L^\dual\otimes\hol K\ar[r]                  &0
			\fullstop
		}
	\end{equation}
	The composition $\beta_2\circ\varphi\circ\alpha_1$ must vanish since it defines a section of $\hol L^{-2}\otimes\hol K$.
	Thereby $\varphi$ induces endomorphisms $\hol L\to\hol L$ and $\hol L^\dual\otimes\hol K\to\hol L^\dual\otimes\hol K$ (the
	vertical dashed morphisms) defining sections of $\hol O$.  That $\varphi$ is an isomorphism shows that both are non-trivial
	and thus non-zero multiples of the identity.  But then the two extensions differ by a non-zero constant in $\Ext^1(\hol
	L^\dual\otimes\hol K,\hol L)$.
\end{proof}

$\roots(S)$ is the set of those line bundles $\hol R\in\Pic^0(S)$ that define a (unique) bundle $\hol A_{\hol R}$ and
$\queer(S)$ is the set of those line bundles $\hol L\in\Pic^0(S)$ that do \emph{not} define a unique bundle $\hol E_{\hol L}$.
In the following we always imply $\hol R\in\roots(S)$ and $\hol L\in\hol\Pic^0(S)\setminus\queer(S)$ when we write $\hol A_{\hol
R}$ and $\hol E_{\hol L}$ respectively.

\begin{remark}
	\label{rem:RQ}
	Using remark~\ref{rem:sections} and evaluating the first \name{Chern} class, it is not difficult to see that the above sets
	have the following form on the different class $\VII_0^1$ surfaces:
	\begin{itemize}
		\item For $S$ the half \name{Inoue} surface, $\roots(S)=\sqrt\hol O=\{\hol O,\hol F\}$ and $\queer(S)=\sqrt\hol F$.
		\item For $S$ an \name{Enoki} or the parabolic \name{Inoue} surface, 
			\begin{equation}
				\label{eq:roots}
				\roots(S)=\{\hol M\otimes\hol O(rC)\colon\hol M\in\sqrt\hol O\cup\sqrt{\hol O(C)},\;r\in\N\}
				\fullstop
			\end{equation}
		\item For $S$ an \name{Enoki} surface, $\queer(S)=\emptyset$
		\item For $S$ the parabolic \name{Inoue} surface, $\queer(S)=\roots(S)\cup\sqrt{\hol O(-C)}$
	\end{itemize}
\end{remark}

In particular, since every line bundle in $\Pic^0(S)$ has exactly two square roots, the above sets are finite or countable so
that the bundles $\hol E_{\hol L}$ with $\hol L\in\Pic^0(S)\setminus\queer(S)$ represent the generic case among the filtrable
bundles of type \eqref{eq:fix}.

\medskip

We now restrict our attention to simple bundles.  Simplicity assures that the resulting moduli space is a complex analytic
space.

\begin{definition}
	A holomorphic vector bundle $\hol E$ is called \dfn{simple} if the only holomorphic endomorphisms of $\hol E$ are multiples
	of the identity.
\end{definition}

\begin{proposition}
	\label{prop:simple}
	\begin{enumerate}
		\item \label{prop:simple:1}
			The central terms of trivial extensions of type \eqref{eq:EA} are never simple.
		\item \label{prop:simple:2}
			The bundles $\hol E_{\hol L}$, $\hol L\in\Pic^0(S)\setminus\queer(S)$, are simple.
		\item \label{prop:simple:3}
			For $\hol L\in\queer(S)$ the central terms of non-trivial extensions \eqref{eq:E} are not simple.
		\item \label{prop:simple:4}
			A bundle $\hol A_{\hol R}$ is simple if $\hol R\in\roots(S)\setminus\queer(S)$.
			\footnote{%
				Later on we will see that this is actually an ``if and only if''.
			}
	\end{enumerate}
	Moreover, every simple filtrable holomorphic bundle of type \eqref{eq:fix} is isomorphic to either a bundle $\hol E_{\hol
	L}$ for some $\hol L\in\Pic^0(S)\setminus\queer(S)$ or to a bundle $\hol A_{\hol R}$ for some $\hol R\in\roots(S)$.
\end{proposition}

\begin{proof}
	(\ref{prop:simple:1}) is evident.  To prove (\ref{prop:simple:2}) and (\ref{prop:simple:3}) regard diagram \eqref{diag:E12}
	for $\hol E_1=\hol E_2\eqcolon\hol E$ and an \emph{endo}morphism $\varphi\colon\hol E\to\hol E$.  As in the proof of
	proposition~\ref{prop:param9n} $\varphi$ induces endomorphisms $\hol L\to\hol L$ and $\hol L^\dual\otimes\hol K\to\hol
	L^\dual\otimes\hol K$ (the vertical dashed morphisms) which must be multiples of the identity since they define sections in
	$\hol O$.  Let the latter one be $\zeta\id_{\hol L^\dual\otimes\hol K}$ with $\zeta\in\C$.  Then we can substitute $\varphi$
	by $\varphi-\zeta\id_{\hol E}$ to obtain the diagram
	\[
		\xymatrix{
			0\ar[r]&\hol L\ar[r]\ar@{-->}[d]&\hol E\ar[r]^{\beta_1\quad}\ar[d]^{\varphi-\zeta\id_{\hol E}}\ar@{..>}[dl]&\hol L^\dual\otimes\hol K\ar[r]\ar@<4pt>@{-->}[d]^0&0\\
			0\ar[r]&\hol L\ar[r]^{\alpha_2} &\hol E\ar[r]                                                              &\hol L^\dual\otimes\hol K\ar[r]                    &0
		}
	\]
	where this time the endmorphism $\hol L^\dual\otimes\hol K\to\hol L^\dual\otimes\hol K$ on the right is zero.  Therefore
	$\varphi-\zeta\id_{\hol E}$ factorises through $\alpha_2$.  But now the endomorphism $\hol L\to\hol L$ on the left must be
	zero too.  Indeed, if not, it would be an isomorphism and its inverse composed with the morphism $\hol E\to\hol L$ would
	define a splitting of the first extension.  This induces yet another morhpism $\sigma\colon\hol L^\dual\otimes\hol K\to\hol
	L$ from the bundle $\hol L^\dual\otimes\hol K$ in the upper extension to the bundle $\hol L$ in the lower extension (not
	indicated).  This morphism defines an element of $H^0(\hol L^2\otimes K^\dual)$ and is zero if and only if
	$\varphi-\zeta\id=\alpha_2\circ\sigma\circ\beta_1$ is.  This demonstrates (\ref{prop:simple:2}) and (\ref{prop:simple:3}).

	The proof of (\ref{prop:simple:4}) is analogue.  In the corresponding diagram
	\[
		\xymatrix{
			0\ar[r]&\hol R^\dual\otimes\hol K\ar^-{\alpha_1}[r]\ar@<4pt>@{-->}[d]&\hol A_{\hol R}\ar[r]\ar@<-4pt>^\varphi[d]&\hol R\ar[r]\ar@{-->}[d]^{\zeta\id_{\hol R}}&0\\
			0\ar[r]&\hol R^\dual\otimes\hol K\ar                              [r]&\hol A_{\hol R}\ar^{\beta_2}[r]           &\hol R\ar[r]                                &0
		}
	\]
	for a bundle endomorphism $\varphi\colon\hol A_{\hol R}\to\hol A_{\hol R}$ the composition
	$\beta_2\circ\varphi\circ\alpha_1$ is zero by hypothesis since it defines a section of $\hol R^2\otimes\hol K^\dual$.  As
	before we can substitute this diagram by
	\[
		\xymatrix{
			0\ar[r]&\hol R^\dual\otimes\hol K\ar[r]\ar@<4pt>@{-->}[d]_0&\hol A_{\hol R}\ar[r]\ar@<-4pt>^{\varphi-\zeta\id}[d]\ar@{..>}[dl]&\hol R\ar[r]\ar@{-->}[d]^0&0\\
			0\ar[r]&\hol R^\dual\otimes\hol K\ar[r]                    &\hol A_{\hol R}\ar[r]                                             &\hol R\ar[r]              &0
			\fullstop
		}
	\]
	Concluding as above we have $\varphi=\zeta\id$ since $H^0(\hol R^{-2}\otimes\hol K)=0$.
	
	The last statement is now a consequence of corollary~\ref{cor:extensions} and proposition~\ref{prop:param9n}.
\end{proof}

To obtain a \emph{bijective} parametrisation of simple filtrable bundles of type \eqref{eq:fix} we will have to determine
possible isomorphisms of the forms
\begin{equation}
	\label{eq:iso-types}
	\hol E_{\hol L'}\isom\hol E_{\hol L} \qquad
	\hol A_{\hol R'}\isom\hol A_{\hol R} \qquad
	\hol A_{\hol R}\isom\hol E_{\hol L}
	\fullstop
\end{equation}
Regarding the defining extensions \eqref{eq:E_L&A_R} and corollary~\ref{cor:extensions}, these are given by holomorphic bundle
embeddings $\hol L'\hookrightarrow\hol E_{\hol L}$, $\hol R'\hookrightarrow\hol A_{\hol R}$ and $\hol L\hookrightarrow\hol
A_{\hol R}$.

A line bundle extension $0\to\hol M\to\hol E\to\hol O\to0$ is determined by the image $\delta_h(1)$ of the constant $1$ section
in $\hol O$ under the connecting operator $\delta_h\colon H^0(\hol O)\to H^1(\hol M)$ in the associated cohomology sequence.
Given, in addition, a divisor $D>0$ on $S$ there is a second connecting operator $\delta_v:H^0(\hol M_D(D))\to H^1(\hol M)$ from
the cohomology sequence associated to the short exact sequence 
\begin{equation}
	\label{eq:M_D(D)}
	0\quad\longrightarrow\quad\hol M\quad\longrightarrow\quad \hol M(D)\quad\longrightarrow\quad\hol M_D(D)\quad\longrightarrow\quad0
	\fullstop
\end{equation}
This sequence is the defining sequence for $\hol M_D(D)$ where we write $\hol M(D)$ for $\hol M\otimes\hol O(D)$ and $\hol M_D$
for the restriction of $\hol M$ to $D$, i.\,e.\ $\hol M_D\coloneq\hol M\otimes\hol O_D$.

In \cite{Teleman05} we find the following criterion:

\begin{proposition}
	\label{prop:linesubbundle}
	With the above notation, the natural map $\hol O(-D)\to\hol O$ can be lifted to a bundle embedding
	\begin{equation}
		\label{eq:bundle_embedding}
		\xymatrix{
			&&\hol O(-D)\ar@{-->}[d]\ar[dr]\\
			0\ar[r]&\hol M\ar[r]&\hol E\ar[r]&\hol O\ar[r]&0
		}
	\end{equation}
	if and only if there exists a section $\sigma\in H^0(\hol M_D(D))$ defining a trivialisation $\hol M_D(D)\isom\hol O_D$ such
	that $\delta_h(1)$=$\delta_v(\sigma)$
\end{proposition}

Applying this criterion to the extensions \eqref{eq:E_L&A_R} yields the following

\begin{corollary}
	\label{cor:isos}
	\begin{enumerate}
		\item \label{cor:isos:E}
			$\hol E_{\hol L}\isom\hol E_{\hol L'}$ if and only if $\hol L'\isom\hol L$.
		\item \label{cor:isos:A}
			Suppose $\hol R'\not\isom\hol R$ and that $\hol A_{\hol R}$ and $\hol A_{\hol R'}$ are simple.  Then $\hol A_{\hol
			R'}\isom\hol A_{\hol R}$ if and only if there exists a divisor $D>0$ with $\hol R\otimes\hol R'\isom\hol K(D)$ and
			$\hol R'_D\isom\hol R_D$.
		\item \label{cor:isos:EA}
			Suppose $\hol A_{\hol R}$ is simple.  Then $\hol A_{\hol R}\isom\hol E_{\hol L}$ if and only if there exists a
			divisor $D>0$ with $\hol L\isom\hol R(-D)$ and $\hol R_D^2\isom\hol K_D(D)$.
	\end{enumerate}
\end{corollary}

\begin{proof}
	We can show the first statement without using proposition~\ref{prop:linesubbundle}.  Take two non-isomorphic bundles $\hol
	L$ and $\hol L'$.  Then either $\hol L^\dual\otimes\hol L'$ or $\hol L\otimes\hol L'^\dual$ has only trivial sections,
	c.\,f.\ remark~\ref{rem:sections}.  We can assume the latter by possibly interchanging $\hol L$ and $\hol L'$.  Let now
	$\varphi\colon\hol E_{\hol L}\to\hol E_{\hol L'}$ be an isomorphism between the corresponding bundles $\hol E_{\hol L}$ and
	$\hol E_{\hol L'}$ and regard the following diagram:
	\[
		\xymatrix{
			0\ar[r]&\hol L \ar^{\alpha_1}[r]&\hol E_{\hol L }\ar[r]\ar^\varphi[d]\ar@{..>}[dl]&\hol L ^\dual\otimes\hol K\ar[r]\ar@{-->}^0[d]&0\\
			0\ar[r]&\hol L'\ar^{\alpha_2}[r]&\hol E_{\hol L'}\ar^{\beta_2\quad}[r]            &\hol L'^\dual\otimes\hol K\ar[r]              &0
			\fullstop
		}
	\]
	The composition $\beta_2\circ\varphi\circ\alpha_1$ vanishes since it defines a section of $\hol L^\dual\otimes\hol
	L'^\dual\otimes\hol K$.  Thus $\varphi$ induces a morphism $\hol L^\dual\otimes\hol K\to\hol L'^\dual\otimes\hol K$ (the
	vertical dashed morphism).  This defines a section of $\hol L\otimes\hol L'^\dual$ which is zero by the above choice of
	$\hol L$ and $\hol L'$.  Consequently $\varphi$ factorises through $\alpha_2$, showing that it can not be an isomorphism.
	This proves the first statement.

	To prove the second statement we can assume that $\hol R^\dual\otimes\hol R'$ does only admit trivial sections by possibly
	interchanging $\hol R$ and $\hol R'$, c.\,f.\ remark~\ref{rem:sections}.  Now observe that an isomorphism $\hol A_{\hol
	R'}\isom\hol A_{\hol R}$ gives, after tensorising the defining extensions for $\hol A_{\hol R}$ and $\hol A_{\hol R'}$ by
	$\hol R^\dual$, a bundle embedding $\alpha$,
	\begin{equation}
		\label{diag:embedding}
		\xymatrix{
			&\quad\hol M\ar@{=}@<1ex>[d]&\hol R^\dual\otimes\hol R'^\dual\otimes\hol K\ar_\alpha[d]\ar@{-->}[dr]& \\
			0\ar[r]&\hol R^{-2}\otimes\hol K\ar[r]&\hol R^\dual\otimes\hol A_{\hol R}\ar[r]&\hol O\ar[r]&0
			\comma
		}
	\end{equation}
	and thus a bundle morphism $\hol R^\dual\otimes\hol R'^\dual\otimes\hol K\to\hol O$.  If it were trivial, $\alpha$ would
	induce a morphism $\hol R^\dual\otimes\hol R'^\dual\otimes\hol K\to\hol R^{-2}\otimes\hol K$ defining a section of $\hol
	R^\dual\otimes\hol R'$ which is zero by assumption.  This would contradict the fact that $\alpha$ is a bundle embedding.  So
	the morphism $\hol R^\dual\otimes\hol R'^\dual\otimes\hol K\to\hol O$ is non-trivial, showing the existence of a divisor
	$D\ge0$ with $\hol R\otimes\hol R'\isom\hol K(D)$.  We have $D\not=0$ because otherwise this would give a splitting of the
	extension defining $\hol A_{\hol R}$, but $\hol A_{\hol R}$ is simple by hypothesis.  Proposition~\ref{prop:linesubbundle}
	applied to $\hol M=\hol R^{-2}\otimes\hol K$ now yields $\hol R^2_D\isom\hol K_D(D)$ or equivalently $\hol R_D\isom\hol
	R'_D$.

	Conversely, suppose $\hol R\otimes\hol R'\isom\hol K(D)$ and $\hol R_D\isom\hol R'_D$.  Again we can assume that $\hol
	R^\dual\otimes\hol R'$ does only admit trivial sections by possibly interchanging $\hol R$ and $\hol R'$.  Consider the
	short exact sequence \eqref{eq:M_D(D)} for $\hol M=\hol R^{-2}\otimes\hol K$ and regard the associated long exact sequence
	\[
		\ldots
		\to H^0(\hol R^\dual\otimes\hol R')
		\to H^0(\hol O_D)
		\stackrel{\delta_v}{\longrightarrow}H^1(\hol R^{-2}\otimes\hol K)
		\to\ldots
	\]
	We have $h^0(\hol R^\dual\otimes\hol R')=0$ by assumption, so the connecting operator $\delta_v$ is injective.  As we saw in
	the proof of proposition~\ref{prop:param9n}, $h^1(\hol R^{-2}\otimes\hol K)=\dim\Ext^1(\hol R,\hol R^\dual\otimes\hol K)=1$.
	Together with $h^0(\hol O_D)\ge1$ this shows that $\delta_v$ is an isomorphism and $h^0(\hol O_D)=1$.  Thus the preimage
	$\sigma$ of $\delta_h(1)\in H^1(\hol R^{-2}\otimes\hol K)$ under $\delta_v$ is a constant section of $\hol M_D(D)\isom\hol
	O_D$ and therefore defines a trivialisation.  Applying proposition~\ref{prop:linesubbundle} to $\hol M=\hol
	R^{-2}\otimes\hol K$ now gives a line bundle inclusion $\alpha$ in \eqref{diag:embedding}.  By
	corollary~\ref{cor:extensions} the resulting bundle embedding $\hol R'^\dual\otimes\hol K\to\hol A_{\hol R}$ extends to an
	extension $0\to\hol R'^\dual\otimes\hol K\to\hol A_{\hol R}\to\hol R'\to0$.  It is non-trivial because $A_{\hol R}$ is
	simple.  Then by the definition of $A_{\hol R'}$ we have $\hol A_{\hol R'}\isom\hol A_{\hol R}$.
	
	The proof of the last statement is analogue because in this case the corresponding diagram is
	\[
		\xymatrix{
			&\quad\hol M\ar@{=}@<1ex>[d]&\hol R^\dual\otimes\hol L\ar^\alpha[d]\ar@{-->}[dr]& \\
			0\ar[r]&\hol R^{-2}\otimes\hol K\ar[r]&\hol R^\dual\otimes\hol A_{\hol R}\ar[r]&\hol O\ar[r]&0
		}
	\]
	and the cohomology sequence of \eqref{eq:M_D(D)} for $\hol M=\hol R^{-2}\otimes\hol K$ reads
	\[
		\ldots
		\to H^0(\hol R^\dual\otimes\hol L^\dual\otimes\hol K)
		\to H^0(\hol O_D)
		\stackrel{\delta_v}{\longrightarrow}H^1(\hol R^{-2}\otimes\hol K)
		\to\ldots
		\fullstop
	\]
	But in this case $\hol R^\dual\otimes\hol L^\dual\otimes\hol K$ does not admit non-trivial sections.
\end{proof}

We will now examine the above criteria on each type of class $\VII_0^1$ surfaces.  For the half \name{Inoue} surface we first
need the following fact.

\begin{lemma}
	\label{lem:K_C(C)}
	A singular rational curve $C$ with one node on a complex surface satisfies $\hol K_C(C)\isom\hol O_C$.
\end{lemma}

\begin{proof}
	Note that $\hol K_C(C)$ is the dualising bundle of $C$ which is independent of the particular embedding of $C$
	\cite{BHPvV} and we can embed $C$ as a cubic in $\C\rm P^2$.  But there $\hol K=\hol O(-3)$ and $\hol O(C)=\hol O(3)$ so
	that $\hol K(C)$ is already trivial.
\end{proof}

\begin{theorem}
	\label{theo:half}
	For $S$ the half \name{Inoue} surface, there is an isomorphism
	\begin{equation}
		\label{eq:iso:half}
		\hol A_{\hol O}\isom\hol A_{\hol F}
	\end{equation}
	and the filtrable simple holomorphic bundles of type \eqref{eq:fix} are bijectively parametrised by
	$\big(\Pic^0(S)\setminus\sqrt{\hol F}\:\big)\amalg\{0\}$, mapping $\hol L\mapsto\hol E_{\hol L}$ and $0\mapsto\hol A_{\hol
	O}\isom\hol A_{\hol F}$.
\end{theorem}

\begin{proof}
	The isomorphism $\hol A_{\hol F}\isom\hol A_{\hol O}$ follows directly from corollary~\ref{cor:isos}(\ref{cor:isos:A}) and
	\eqref{eq:half} together with the lemma.  Note that by remark~\ref{rem:RQ} there are no further bundles of the form $\hol
	A_{\hol R}$.
	
	The bundles $\hol E_{\hol L}$ are simple by proposition~\ref{prop:simple} as well as is $\hol A_\hol O$ because $\hol
	O\not\in\queer(S)=\sqrt{\hol F}$.  By corollary~\ref{cor:isos} the bundles $\hol E_{\hol L}$ are pairwise non-isomorphic and
	there can be no isomorphism $\hol A_{\hol O}\isom\hol E_{\hol L}$.  Indeed, taking the first \name{Chern} class of $\hol
	L\isom\hol R(-D)$ shows $c_1(\hol O(D))=0$, contradicting $D\not=0$.  This shows injectivity.  Surjectivity follows from
	proposition~\ref{prop:simple}.
\end{proof}

\begin{theorem}
	\label{theo:parabolic}
	For $S$ the parabolic \name{Inoue} surface there are isomorphisms 
	\begin{equation}
		\label{eq:iso:parabolic}
		\hol A_{\hol R}\isom\hol R(-E)\oplus\hol R^\dual(-C)
		\qquad
		\hol R\in\roots(S)
		\comma
	\end{equation}
	so the $\hol A_{\hol R}$ are not simple.  The filtrable simple bundles of type \eqref{eq:fix} are bijectively parametrised
	by $\Pic^0(S)\setminus\queer(S)$, mapping $\hol L\mapsto\hol E_{\hol L}$.
\end{theorem}

\begin{proof}
	To show \eqref{eq:iso:parabolic}, take a bundle $\hol R\in\roots(S)$ with $\hol R^2\isom\hol O(rC)$ for some $r\in\N$.
	Using $\hol K\isom\hol O(-E-C)$ we get
	\[
		(\hol K\otimes\hol R^\dual)^\dual\;\otimes\;\big(\hol R(-E)\oplus\hol R^\dual(-C)\big)
		\;\;=\;\;
		\hol O\big((r+1)C\big)\;\oplus\;\hol O(E)
		\fullstop
	\]
	Since $C\cap E=\emptyset$, this bundle admits a non-vanishing section giving rise to a bundle embedding $\hol K\otimes\hol
	R^\dual\hookrightarrow\hol R(-E)\oplus\hol R^\dual(-C)$.  But, as one easily checks, $\hol R(-E)\oplus\hol R^\dual(-C)$ is
	of type \eqref{eq:fix}.  So by corollary~\ref{cor:extensions} this inclusion extends to
	\[
		0\quad\longrightarrow\quad\hol K\otimes\hol R^\dual\quad\longrightarrow\quad\hol R(-E)\oplus\hol R^\dual(-C)\quad\longrightarrow\quad\hol R\quad\longrightarrow\quad0
		\fullstop
	\]
	Assume this extension splits, i.\,e.\ $\hol R(-E)\oplus\hol R^\dual(-C)\isom(\hol K\otimes\hol R^\dual)\oplus\hol R$.
	Tensorising with $\hol R^\dual$ gives $\hol O(-E)\oplus\hol O(-(r+1)C)\isom(\hol K\otimes\hol R^{-2})\oplus\hol O$ which is
	impossible because the left hand side admits no non-trivial sections while the right hand side does.  Therefore the above
	extension is non-trivial and determines the isomorphism \eqref{eq:iso:parabolic} by the very definition of $\hol A_{\hol
	R}$.  The rest follows from \ref{prop:simple} and \ref{cor:isos}(\ref{cor:isos:E}).
\end{proof}

To apply corollary~\ref{cor:isos} in the remaining case of an \name{Enoki} surface, we need the following generalisation of
lemma~\ref{lem:K_C(C)} for \name{Enoki} surfaces.

\begin{lemma}
	On an \name{Enoki} surface one has $\hol K_{rC}(C)\isom\hol O_{rC}$ for $r\in\N\setminus\{0\}$.
\end{lemma}

\begin{proof}
	We prove by induction on $r$.  For $r=1$ this is just lemma~\ref{lem:K_C(C)}, so let us suppose $\hol K_{rC}(C)\isom\hol
	O_{rC}$ for some $r\ge1$.  Restricting a holomorphic line bundle $\hol M$ from $(r+1)C$ to $C$ gives the following exact
	sequence \cite{BHPvV}:
	\begin{equation}
		\label{eq:restriction}
		0\quad\longrightarrow\quad\hol M_{rC}(-C)\quad\longrightarrow\quad \hol M_{(r+1)C}\quad\longrightarrow\quad\hol M_C\quad\longrightarrow\quad0
		\fullstop
	\end{equation}
	Considering this sequence for $\hol M=\hol O$ and $\hol M=\hol K(C)$ respectively and taking into account
	lemma~\ref{lem:K_C(C)} as well as the induction hypothesis gives the following two extensions of $\hol O_{rC}(-C)$ by $\hol
	O_C$:
	\[
		\xymatrix{
			0\ar[r]&\hol O_{rC}(-C)\ar[r]\ar@{=}[d]&\hol O_{(r+1)C}   \ar[r]&\hol O_C   \ar[r]\ar@{=}[d]&0\\
			0\ar[r]&\hol K_{rC}    \ar[r]          &\hol K_{(r+1)C}(C)\ar[r]&\hol K_C(C)\ar[r]          &0
			\fullstop
		}
	\]
	But the set of isomorphism classes of holomorphic line bundles on $(r+1)C$ extending $\hol O_{rC}(-C)$ by $\hol O_C$ can be
	identified with $H^1(\hol O_C(-rC))$ \cite{Drezet}%
	\footnote{%
		The proof in \cite{Drezet} goes through for non-algebraic surfaces as well.
	}.
	So all we have to verify is $H^1(\hol O_C(-rC))=0$.  To this aim consider the part
	\[
		\ldots\to H^1\big(\hol O(-rC)\big)\to H^1\big(\hol O_C(-rC)\big)\to H^2\big(\hol O(-(r+1)C)\big)\to\ldots
	\]
	of the exact cohomology sequence associated to the short exact sequence \eqref{eq:M_D(D)} for $\hol M=\hol O(-(r+1)C)$ and
	$D=C$.  By \name{Serre} duality we have $H^2(\hol O(-(r+1)C))\isom H^0(\hol K((r+1)C))=0$ by remark~\ref{rem:sections} since
	$c_1(\hol O(C))=0$ on an \name{Enoki} surface.  Also $H^0(\hol K(rC))=0$ and formula \eqref{eq:RRSD} for $\hol M=\hol
	O(-rC)$ shows $H^1(\hol O(-rC))=0$.  This proves $H^1(\hol O_C(-rC))=0$ and therefore $\hol K_{(r+1)C}(C)\isom\hol
	O_{(r+1)C}$.
\end{proof}

\begin{theorem}
	\label{theo:Enoki}
	On an \name{Enoki} surface $S$ we have isomorphisms%
	\footnote{
		$\hol E_{\hol R^\dual(-C)}$ is well defined since $\hol R^\dual(-C)\in\Pic^0(S)\setminus\queer(S)$, as one easily
		verifies.
	}
	\begin{equation}
		\label{eq:iso:Enoki}
		\hol A_{\hol R}\isom\hol E_{\hol R^\dual(-C)}
		\qquad
		\hol R\in\roots(S)
	\end{equation}
	and the filtrable simple bundles of type \eqref{eq:fix} are bijectively parametrised by $\Pic^0(S)$, mapping $\hol
	L\mapsto\hol E_{\hol L}$.
\end{theorem}

\begin{proof}
	First recall that divisors on an \name{Enoki} surface are multiples of the curve $C$.  Therefore $\hol R^2\isom\hol O(rC)$
	for some $r\in\N$ if $\hol R\in\roots(S)$ and the above lemma shows the existence of a divisor $D=(r+1)C$ with $\hol
	K_D(D)\isom\hol O_D(rC)\isom\hol R^2_D$.  Corollary~\ref{cor:isos}(\ref{cor:isos:EA}) thus gives an isomorphism $\hol
	A_{\hol R}\isom\hol E_{\hol L}$ with $\hol L\isom\hol R\big(-(r+1)C\big)=\hol R^\dual(-C)$.  Remarking that
	$\queer(S)=\emptyset$ for \name{Enoki} surfaces, the rest follows from \ref{prop:simple} and
	\ref{cor:isos}(\ref{cor:isos:E}). 
\end{proof}

Resuming, we saw that with exception of the bundle $\hol A_{\hol O}\isom\hol A_{\hol F}$ on the half \name{Inoue} surface, every
filtrable simple holomorphic bundle of type \eqref{eq:fix} on a class $\VII_0^1$ surface $S$ is of the form $\hol E_{\hol L}$
with $\hol L\in\Pic^0(S)\setminus\queer(S)$.  Taking into account remark~\ref{rem:RQ}, a bundle of the form $\hol A_{\hol R}$ is
simple if and only if $\hol R\in\roots(S)\setminus\queer(S)$.

\section{The local structure of the moduli space}
\label{sec:local}

\begin{definition}
	We denote by
	\[
		\mathcal M^{\rm s}(S)
			\coloneq\{\text{$\hol E$ simple holomorphic structure on $E$}\colon\det\hol E\isom\hol K\}/\Gamma\big(S,\GL(E)\big)
	\]
	the \dfn{moduli space of simple holomorphic bundles} of type \eqref{eq:fix} on $S$.
\end{definition}

This is a (possibly non-\name{Hausdorff}) complex space.  The local structure of this moduli space is given by the following
proposition whose proof is a straightforward generalisation of the case $\roots(S)=\sqrt{\hol O}$ in \cite{Teleman}.

\begin{proposition}
	\label{prop:local}
	\begin{enumerate}
		\item
			\label{prop:local:1}
			If $\hol L\in\Pic^0(S)\setminus\big(\roots(S)\cup\queer(S)\big)$ then $\mathcal M^{\rm s}(S)$ is a smooth complex
			curve $C_{\hol L}$ in a neighbourhood of $\hol E_{\hol L}$ given by $\hol L'\mapsto\hol E_{\hol L'}$.
		\item
			\label{prop:local:2}
			If $\hol R\in\roots(S)\setminus\queer(S)$ then $\mathcal M^{\rm s}(S)$ is the intersection of two complex curves
			$C_{\hol R}$ and $C'_{\hol R}$ in a neighbourhood of $\hol E_{\hol R}$ where $C_{\hol R}$ is given by $\hol
			L'\mapsto\hol E_{\hol L'}$.
		\item
			\label{prop:local:3}
			If $\hol R\in\roots(S)\setminus\queer(S)$ then $\mathcal M^{\rm s}(S)$ is a smooth complex curve $C''_{\hol R}$ in
			a neighbourhood of $\hol A_{\hol R}$.
		\item
			\label{prop:local:4}
			The points $\hol E_{\hol R}$ and $\hol A_{\hol R}$ are not separable.  More precisely, we find neighbourhoods $U'$,
			$U''$ of $\hol E_{\hol R}$ and $\hol A_{\hol R}$ respectively with $\big(C'_{\hol R}\setminus\{\hol E_{\hol
			R}\}\big)\cap U'=\big(C''_{\hol R}\setminus\{\hol A_{\hol R}\}\big)\cap U''$.
		\item
			\label{prop:local:5}
			$\mathcal M^{\rm s}(S)$ is a smooth complex curve in a neighbourhood of every non-filtrable bundle.
	\end{enumerate}
	Moreover, $\mathcal M^{\rm s}(S)$ is regular in all smooth points.
\end{proposition}

\begin{figure}[h]
	\centering
	\includegraphics{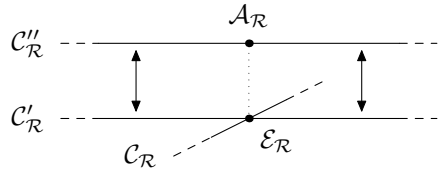}
	\caption{Local structure of the moduli space at $\hol E_{\hol R}$ and $\hol A_{\hol R}$}
	\label{fig:local}
\end{figure}

We have depicted this situation in figure~\ref{fig:local} (dividing real dimensions by two).  The vertical arrows symbolise
identification of the corresponding curves with exception of the points joined by the dotted line.  We can regard the curves
$C'_{\hol R}$ and $C''_{\hol R}$ as one single curve with a double point consisting of $\hol E_{\hol R}$ and $\hol A_{\hol R}$.
This ``curve'' is smooth at the point $\hol A_{\hol R}$ but is transversely crossed by the curve $C_{\hol R}$ at the point $\hol
E_{\hol R}$.

In the case of the parabolic \name{Inoue} surface and an \name{Enoki} surface the above theorem determines completely the
structure of the moduli space in a neighbourhood of every filtrable bundle.  Recall that for the parabolic \name{Inoue} surface
$\roots(S)\setminus\queer(S)=\emptyset$ so that the situation is particularly simple:  Theorem~\ref{theo:parabolic} actually
establishes an \emph{isomorphism} between $\Pic^0(S)\setminus\queer(S)$ and the filtrable part of the moduli space of simple
bundles, given by $\hol L\mapsto\hol E_{\hol L}$.  For an \name{Enoki} surface the isomorphisms $\hol A_{\hol R}\isom\hol
E_{\hol R^\dual(-C)}$ immediately tell us that $C''_{\hol R}=C_{\hol R^\dual(-C)}$.  In the remaining case of the half
\name{Inoue} surface the situation is slightly more complicated.  We can not yet identify the curves $C'_{\hol R}$ and
$C''_{\hol R}$ and will do this in section~\ref{sec:centres}.

\section{Stability}
\label{sec:stable}

The moduli space important for gauge theory is the moduli space of stable holomorphic bundles and is a \emph{\name{Hausdorff}}
complex space.  Stability is defined with respect to a \dfn{\name{Gauduchon} metric} $g$ on $S$ which is a \name{Hermit}ian
metric whose associated $(1,1)$-form $\omega_g$ verifies $\partial\bar\partial\omega_g=0$.  Such a metric always exists and
allows one to define the \dfn{degree map} by
\[
	\begin{array}{rcl}
		\deg\colon\Pic(S)&\longrightarrow&\mathbb R\\
		\hol L&\longmapsto&\deg\hol L\coloneq{\displaystyle\int\limits_S}c_1(\hol L,A_h)\wedge\omega_g
		\comma
	\end{array}
\]
where $c_1(\hol L,A_h)$ is the first \name{Chern} form associated to the \name{Chern} connection $A_h$ of a \name{Hermit}ian
metric $h$ in $\hol L$ (i.\,e.\ locally $c_1(\hol L,A_h)=\partial\bar\partial\log h$).  This is a \name{Lie} group morphism and
independent of the particular choice of $h$.  Note that on non-\name{K\"ahler} surfaces the degree map is never a topological
invariant and therefore non-constant on $\Pic^0(S)$ \cite{Luebke&Teleman}.

\begin{examples}
	\label{ex:degree}
	\begin{enumerate}
		\item
			\label{ex:degree:F}
			$\deg\hol O=\deg\hol F=0$ because the square roots of $\hol O$ are torsion elements in $\Pic^0(S)$.
		\item
			\label{ex:degree:D}
			$\deg\hol O(D)=\vol D>0$ for any divisor $D>0$ on $S$.  This is a consequence of the
			\name{Poincar\'e}-\name{Lelong} formula \cite{Griffith&Harris}.
		\item
			\label{ex:degree:K}
			On the half and the parabolic \name{Inoue} surface $\deg\hol K<0$ for any \name{Gauduchon} metric.  This follows
			from the previous examples together with \eqref{eq:half} and \eqref{eq:parabolic} respectively.  For \name{Enoki}
			surfaces $\deg\hol K$ attains every value in $\R$ when $g$ varies in the space of \name{Gauduchon} metrics.
			This was shown in \cite{Teleman05}, based on results of \cite{Buchdahl00}.
		\item
			\label{ex:degree:R}
			$\deg\hol R\ge0$ if $\hol R\in\roots(S)$, because $\hol R^2\isom\hol O(D)$ for some divisor $D\ge0$.
		\item
			\label{ex:degree:L}
			Likewise, $\deg\hol L\ge\frac12\deg\hol K$ if $\hol L\in\queer(S)$.
	\end{enumerate}
\end{examples}

Now the (slope-)stability is defined using the \dfn{$g$-slope} of a coherent sheaf $\sheaf S$
\[
	\mu_g(\sheaf S)\coloneq\frac{\deg\det\sheaf S}{\operatorname{rank}\sheaf S}
	\fullstop
\]

\begin{definition}
	\label{def:stable}
	A holomorphic rank two vector bundle $\hol E$ over a complex surface $S$ is called \dfn{$g$-stable} if for every rank one
	subsheaf $\sheaf S\subset\hol E$ we have $\mu_g(\sheaf S)<\mu_g(\hol E)$.
\end{definition}

This definition simplifies in our case to:

\begin{proposition}
	A holomorphic vector bundle $\hol E$ of type \eqref{eq:fix} on a class $\VII_0^1$ surface $S$ is $g$-stable if and only if
	for every holomorphic line-subbundle $\hol L\subset\hol E$ we have $\deg\hol L<\frac12\deg\hol K$.
\end{proposition}

The proof is standard, see for example \cite{Kobayashi}.


Non-filtrable bundles are stable by definition and stable bundles are simple \cite{Kobayashi}, so it remains to examine stability
for simple filtrable bundles (c.\,f.\ section~\ref{sec:filtrable}).

\begin{proposition}
	\label{prop:subbundles}
	\begin{enumerate}
		\item
			\label{prop:subbundles:half}
			On the half \name{Inoue} surface, the bundle $\hol A_{\hol O}\isom\hol A_{\hol F}$ has exactly two holomorphic line
			subbundles, namely $\hol K$ and $\hol F\otimes\hol K$.
		\item
			\label{prop:subbundles:Enoki}
			On an \name{Enoki} surface, the bundles $\hol A_{\hol R}\isom\hol E_{\hol R^\dual(-C)}$ have exactly two holomorphic
			line subbundles, namely $\hol R^\dual\otimes\hol K$ and $\hol R^\dual(-C)$.
		\item
			On an arbitrary class $\VII_0^1$ surface, a bundle $\hol E_{\hol L}$ has no holomorphic line subbundle other than
			$\hol L$ if it does not belong to case~(\ref{prop:subbundles:Enoki}).
	\end{enumerate}
\end{proposition}

\begin{proof}
	By definition, the bundles $\hol E_{\hol L}$ and $\hol A_{\hol R}$ have as holomorphic line subbundles $\hol L$ and
	$\hol R^\dual\otimes\hol K$ respectively.  By corollary~\ref{cor:extensions} every other inclusion of a holomorphic
	line bundle into $\hol E_{\hol L}$ or $\hol A_{\hol R}$ extends to an extension of type \eqref{eq:EA}.  This extension is
	non-trivial if the bundle $\hol E_{\hol L}$ respectively $\hol A_{\hol R}$ is simple and thus determines an isomorphism
	\eqref{eq:iso-types} of the corresponding central terms.  The proposition follows now from the classification
	\eqref{eq:iso:half}--\eqref{eq:iso:Enoki} of all possible such isomorphisms.
\end{proof}

\begin{corollary}
	\label{cor:stability}
	\begin{enumerate}
		\item \label{cor:stability:1}
			On the half \name{Inoue} surface, the bundle $\hol A_{\hol O}\isom\hol A_{\hol F}$ is $g$-stable for any
			\name{Gauduchon} metric $g$.
		\item \label{cor:stability:2}
			On an \name{Enoki} surface, $\hol A_{\hol R}\isom\hol E_{\hol R^\dual(-C)}$ is $g$-stable if and only if
			\[
				\left\{
					\begin{aligned}
						\quad\deg\hol R^\dual(-C)&<\tfrac12\deg K&\text{in case}\quad\deg\hol K&<0\\
						           \tfrac12\deg K&<\deg\hol R    &\text{in case}\quad\deg\hol K&\ge0
					\end{aligned}
				\right.
				\fullstop
			\]
			In either case, one inequality implies the other.
		\item \label{cor:stability:3}
			On an arbitrary class $\VII_0^1$ surface, a bundle $\hol E_{\hol L}$ not belonging to case~(\ref{cor:stability:2})
			is $g$-stable if and only if $\deg\hol L<\tfrac12\deg\hol K$.
	\end{enumerate}
\end{corollary}

\begin{proof}
	Combine the previous proposition with the examples~\ref{ex:degree}.
\end{proof}

\begin{remark}
	We see that always at least one of the two non-separable bundles $\hol E_{\hol R}$ and $\hol A_{\hol R}$ is unstable as it
	should be, for the moduli space of stable bundles is \name{Hausdorff}.
\end{remark}

The degree homomorphism is non-constant on $\Pic^0(S)\isom\C^*$, so the degree corresponds to a non-zero multiple of the
logarithm of the radius in $\C$.  Regarding \ref{cor:stability}(\ref{cor:stability:3}) we fix an isomorphism
$\Pic^0(S)\isom\C^*$ that identifies 
\[
	\Pic^0_{<\varrho}(S)\coloneq\{\hol L\in\Pic^0(S)\colon\deg\hol L<\varrho\}
	\qquad
	\varrho\coloneq\tfrac12\deg\hol K
\]
to a punctured open disc in $\C$ with center $0$, corresponding to $\deg\hol L\to-\infty$.  In view of
\ref{cor:stability}(\ref{cor:stability:2}) we also define the set
\begin{equation}
	\label{eq:unstables}
	\unstable(S)\coloneq\{\hol R^\dual(-C)\in\Pic^0(S)\colon\hol R\in\roots(S),\;\deg\hol R\le\tfrac12\deg K\}
	\fullstop
\end{equation}
From \ref{ex:degree}(\ref{ex:degree:R}) we see that $\unstable(S)=\emptyset$ if $\deg\hol K<0$ --- in particular if $S$ is the
half or the parabolic \name{Inoue} surface.  If $S$ is an \name{Enoki} surface then $\unstable(S)$ is the finite set consisting
of those line bundles $\hol L\in\Pic^0(S)$ with $\deg\hol L<\varrho$ that define an unstable bundle $\hol E_{\hol L}$.  Note
that under the map $\hol R\mapsto\hol R^\dual(-C)$ the set $\unstable(S)$ is in bijection to the set
$\roots_{\le\varrho}(S)\coloneq\roots(S)\cap\Pic^0_{\le\varrho}(S)$ defining singular semistable points $\hol E_{\hol R}$ in the
moduli space.

\begin{corollary}
	\label{cor:stable}
	The filtrable part of the moduli space of $g$-stable holomorphic bundles is bijectively parametrised by
	\begin{itemize}
		\item $\Pic^0_{<\varrho}(S)$ if $S$ is the parabolic \name{Inoue} surface ($\unstable(S)=\emptyset$),
		\item $\Pic^0_{<\varrho}(S)\amalg\{0\}$ if $S$ is the half \name{Inoue} surface ($\unstable(S)=\emptyset$) and
		\item $\Pic^0_{<\varrho}(S)\setminus\:\unstable(S)$ if $S$ is an \name{Enoki} surface,
	\end{itemize}
	mapping $\Pic^0_{<\varrho}(S)\ni\hol L\mapsto\hol E_{\hol L}$ and $0\mapsto\hol A_{\hol O}$.
\end{corollary}

\begin{proof}
	This follows from the above corollary together with the theorems~\ref{theo:half}, \ref{theo:parabolic}, \ref{theo:Enoki} and
	the observation from example~\ref{ex:degree}(\ref{ex:degree:L}) that $\queer(S)\cap\Pic^0_{<\varrho}(S)=\emptyset$.
\end{proof}

\section{The boundary of the moduli space of polystable bundles}
\label{sec:boundary}

We want to compute the moduli spaces of polystable holomorphic bundles of type \eqref{eq:fix} for any class $\VII_0^1$ surface
$S$.  Throughout this section we fix $S$ and omit it in notations.

\begin{definition}
	A holomorphic rank two vector bundle $\hol E$ is \dfn{$g$-polystable} if it is $g$-stable (definition~\ref{def:stable}) or
	if 
	\begin{equation}
		\label{eq:split_polystable}
		\hol E=\hol L\oplus\hol M
		\qquad
		\text{with}
		\quad
		\deg\hol M=\deg\hol L
		\fullstop
	\end{equation}
	In the latter case we call $\hol E$ a \dfn{split $g$-polystable bundle}.  We denote by
	\[
		\mathcal M^{\rm (p)st}\coloneq
			\{\text{$\hol E$ (poly)stable hol.\ str.\ on $E$}\colon\det\hol E\isom\hol K\}
			/\Gamma\big(S,\operatorname{GL}(E)\big)
	\]
	the \dfn{moduli space of (poly)stable holomorphic bundles of type \eqref{eq:fix}}.
\end{definition}

In the previous sections we showed that there is an injection of $\Pic^0_{<\varrho}\setminus\:\unstable$ into the filtrable part
of the moduli space $\mathcal M^{\rm st}$ of stable bundles given by $\hol L\mapsto\hol E_{\hol L}$ (corollary
\ref{cor:stable}) which is holomorphic on $\Pic^0_{<\varrho}\setminus\big(\unstable\cup\roots_{\le\varrho}\big)$ (proposition
\ref{prop:local}).  Now define the closed punctured disc
\[
	\Pic^0_{\le\varrho}\coloneq\{\hol L\in\Pic^0\colon\deg\hol L\le\varrho\}\subset\Pic^0\isom\C^*
	\qquad
	\varrho=\tfrac12\deg\hol K
	\fullstop
\]
Its boundary is the circle $\Pic^0_{=\varrho}(S)$ of line bundles $\hol L\in\Pic^0(S)$ with $\deg\hol L=\tfrac12\deg\hol K$ and
can be mapped to the split polystable bundles by $\hol L\mapsto\hol L\oplus(\hol L^\dual\otimes\hol K)$.

In the following we use results about the gauge theoretical counterpart of our complex geometric moduli space.  Equip the bundle
$E$ with a \name{Hermit}ian metric $h$ and fix a $(\det h)$-unitary connection $a$ in the determinant line bundle $\det E=K$.
Denote by
\[
	\mathcal M^{\rm ASD}\coloneq\{\text{$A$ $h$-unitary connection on $E$}\colon F_A^+=0,\det A=a\}/\Gamma\big(S,\SU(E)\big)
\]
the moduli space of oriented anti-self-dual (ASD) connections.  A connection $A$ is called \dfn{reducible} if there is an
$A$-parallel splitting of $E$ into two line bundles, i.\,e.\ $E=L\oplus M$ and $A=A_L\oplus A_M$ where $A_L$ and $A_M$ are
connections on the line bundles $L$ and $M$ respectively.  We write $(\mathcal M^{\rm ASD})^*$ for the irreducible part of
$\mathcal M^{\rm ASD}$ which is naturally a real analytic space.

The relation between this gauge theoretical moduli space of ASD connections and the complex geometric moduli space of
holomorphic bundles is given by the \name{Kobayashi}-\name{Hitchin} correspondence \cite{Luebke&Teleman}, a natural real
analytic isomorphism
\begin{equation}
	\label{eq:KH}
	{\rm KH}\colon\big(\mathcal M^{\rm ASD}\big)^*\stackrel{\isom}{\longrightarrow}\mathcal M^{\rm st}
\end{equation}
given by mapping the gauge equivalence class $[A]$ of an ASD connection $A$ to the holomorphic structure in $E$ determined by
the corresponding $\bar\partial$-operator $\bar\partial_A$.

Now the second reason for our particular choice of the \name{Chern} classes \eqref{eq:fix_E} of $E$ becomes apparent.  The
moduli space $\mathcal M^{\rm ASD}$ has a natural compactification --- the \name{Uhlenbeck} compactification
\cite{Donaldson&Kronheimer}, constructed by attaching further strata involving moduli spaces $\mathcal M^{\rm ASD}(E_k)$ of
oriented ASD connections on rank two bundles $E_k$ with
\[
	c_1(E_k)=c_1(E)
	\quad\text{and}\quad
	c_2(E_k)=c_2(E)-k
	\comma\qquad
	k=1,2,\ldots
	\fullstop
\]
But in our case \eqref{eq:fix_E} assures that $4c_2(E_k)-c_1(E_k)^2<0$, condition under which the expected dimension
\eqref{eq:dim} of $\mathcal M^{\rm ASD}(E_k)$ is negative and the attached strata in the \name{Uhlenbeck} compactification of
$\mathcal M^{\rm ASD}$ are all empty.  This means that $\mathcal M^{\rm ASD}$ \emph{is already compact} and the irreducible part
$(\mathcal M^{\rm ASD})^*$ of $\mathcal M^{\rm ASD}$ can be compactified by adding only the reducible part.  The latter can be
shown to be the circle $iH^1(S,\R)/2\pi iH^1(S,\Z)$.  In fact, applying the \name{Kobayashi}-\name{Hitchin}-correspondence for
line bundles separately to the line bundles in the splitting \eqref{eq:split_polystable} of a split polystable bundle maps the
circle of split polystable bundles to this circle of reducible connections.

Putting together the above, we get the following commutative diagram
\[
	\xymatrix{
		\Pic^0_{  <\varrho}\setminus\:\unstable\ar[r]\ar[d]&\mathcal M^{\rm  st}\ar^-{\isom}_-{\rm KH}[r]\ar[d]&\Big(\mathcal M^{\rm ASD}\Big)^*\ar[d]\\
		\Pic^0_{\le\varrho}\setminus\:\unstable\ar[r]      &\mathcal M^{\rm pst}\ar^-{\isom}          [r]      &     \mathcal M^{\rm ASD}
	}
	\fullstop
\]
where the vertical arrows are natural inclusions.  Remark that a priori there is no natural topology on the moduli space of
polystable bundles and the bijection $\mathcal M^{\rm pst}\rightarrow\mathcal M^{\rm ASD}$ is only set theoretical.  It is
turned tautologically into a homeomorphism by equipping $\mathcal M^{\rm pst}$ with the induced topology.

\begin{proposition}
	\label{prop:boundary}
	The above inclusion
	$\Pic^0_{\le\varrho}\setminus\:\unstable\hookrightarrow\mathcal M^{\rm pst}$ maps
	$\Pic^0_{\le\varrho}\setminus\big(\unstable\cup\roots_{\le\varrho}\big)$ homeomorphically to an open subspace of $\mathcal
	M^{\rm pst}$.  In particular, if there is no bundle $\hol R\in\roots$ with $\deg\hol R=\varrho$, $\mathcal M^{\rm pst}$
	possesses the structure of a real two-dimensional manifold with boundary in the neighbourhood of the image of the circle
	$\Pic^0_{=\varrho}$.
\end{proposition}

\begin{proof}
	Using the following lemma, we can apply the proof of \cite[prop.~4.4]{Teleman}.  Remark that
	$\Pic^0_{=\varrho}\cap\:\unstable=\emptyset$.
\end{proof}

\begin{lemma}
	Let $\hol E$ be a stable holomorphic bundle of type \eqref{eq:fix} and $\varepsilon>0$ be sufficiently small.  Then a line
	bundle $\hol M\in\Pic^0$ with $H^0(\hol M^\dual\otimes\hol E)\not=0$ and $\varrho-\varepsilon\le\deg\hol M\le\varrho$ is
	unique.
\end{lemma}

\begin{proof}
	The existence of such a line bundle $\hol M$ implies that $\hol E$ is filtrable
	and as in the proof of the equivalence in definition~\ref{def:filtrable} we can construct a non-trivial sheaf morphism $\hol
	M\to\hol L$ to a line subbundle $\hol L$ of $\hol E$.
	So $\hol M\isom\hol L(-D)$ for some divisor $D\ge0$ on $S$.  Since $\hol E$ is stable we
	have $\deg\hol L<\varrho$ and $\vol D=\deg\hol L-\deg\hol M<\varepsilon$.  If we choose $\varepsilon>0$ less than the volume
	of any curve on the surface then $D=0$ and $\hol M\isom\hol L$.  But proposition~\ref{prop:subbundles} shows that a line
	subbundle $\hol L\in\Pic^0$ of $\hol E$ is unique.
\end{proof}

The proof of proposition~\ref{prop:boundary} fails at points $\hol R\in\roots$ with $\deg\hol R=\varrho$, c.\,f.\ \cite[lemma
4.3]{Teleman}.  This can only occur on \name{Enoki} surfaces for \name{Gauduchon} metrics with $\deg\hol K>0$ and we will
account for this situation in the last section when we discuss the structure of the entire moduli space.

\section{Non-filtrable holomorphic bundles}
\label{sec:centres}

The next proposition says that the structure of the moduli space around the origin is the natural one given by the closure
$\Pic^0_{\le\varrho}\cup\:\{0\}$ of $\Pic^0_{\le\varrho}$ in $\C$.

\begin{proposition}[\cite{Teleman}, prop 4.5]
	\label{prop:centre}
	The inclusion $\Pic^0_{\le\varrho}\setminus\:\unstable\hookrightarrow\mathcal M^{\rm pst}$ extends to an inclusion
	\[
		\big(\Pic^0_{\le\varrho}\cup\{0\}\big)\setminus\unstable\hookrightarrow\mathcal M^{\rm pst}
		\comma
	\]
	holomorphic at the centre $0$.  Moreover, $0$ is mapped to a bundle $\hol E$ verifying
	\begin{equation}
		\label{eq:invariance}
		\hol E\otimes\hol F\isom\hol E
		\comma
	\end{equation}
	where $\hol F$ is the (unique) non-trivial square-root of $\hol O$.
\end{proposition}

The invariance property \eqref{eq:invariance} follows from the following lemma in the limit $\hol L\to0$, i.\,e.\ $\deg\hol
L\to-\infty$, since $\deg\hol L\otimes\hol F=\deg\hol L$.

\begin{lemma}
	\label{lem:invariance}
	$\hol E_{\hol L}\otimes\hol F\isom\hol E_{\hol L\otimes\hol F}$ and $\hol A_{\hol R}\otimes\hol F\isom\hol A_{\hol
	R\otimes\hol F}$. 
\end{lemma}

\begin{proof}
	First note that $\hol E_{\hol L}\otimes\hol F$ and $\hol A_{\hol R}\otimes\hol F$ are of type \eqref{eq:fix}.  Tensorise the
	defining extensions for $\hol E_{\hol L}$ and $\hol A_{\hol R}$ by $\hol F$ and compare with the defining extensions for
	$\hol E_{\hol L\otimes\hol F}$ and $\hol A_{\hol R\otimes\hol F}$ respectively.
\end{proof}

This also makes explicit the $\Z_2$ symmetry of the moduli spaces of bundles of type \eqref{eq:fix} under tensorising with the
square roots of $\hol O$.  We see that \eqref{eq:invariance} holds for $\hol A_{\hol O}\isom\hol A_{\hol F}$.

\begin{corollary}
	On the half \name{Inoue} surface, $\hol E$ is the filtrable bundle $\hol A_{\hol O}$.  On an \name{Enoki} or the parabolic
	\name{Inoue} surface $\hol E$ is a non-filtrable bundle.
\end{corollary}

\begin{proof}
	Suppose $\hol E_{\hol L}\otimes\hol F\isom\hol E_{\hol L}$ for $S$ an arbitrary class $\VII_0^1$ surface.  Then $\hol
	E_{\hol L\otimes\hol F}\isom\hol E_{\hol L}$ by lemma~\ref{lem:invariance} and thus $\hol L\otimes\hol F\isom\hol L$ by
	corollary~\ref{cor:isos}, contradicting the non-triviality of $\hol F$.  Therefore, either $\hol E$ is non-filtrable or $S$
	is the half \name{Inoue} surface and $\hol E\isom\hol A_{\hol O}$.  $\hol E$ cannot be non-filtrable on the half
	\name{Inoue} surface because this would imply that $\hol A_{\hol O}$ lies on another component of the moduli space.  This is
	excluded by corollary~\ref{cor:connected} below.
\end{proof}

\begin{remark}
	One can show that \eqref{eq:invariance} implies that the pull-back of $\hol E$ to a double cover of $S$ splits into a sum of
	two line bundles.
\end{remark}

For a complete description it only remains to show that our moduli spaces do not contain further connected components.
Non-filtrable bundles are stable by definition and we saw that all unstable filtrable bundles lie on the component we already
described.  Thus another component would be contained in the moduli space of polystable bundles and therefore be compact.  But
\name{M.~Toma} showed that this is impossible on blown-up primary \name{Hopf} surfaces \cite{Toma06} and we know that every
class $\VII_0$ surface containing a global spherical shell --- in particular every class $\VII_0^1$ surface --- is a
degeneration of blown-up primary \name{Hopf} surfaces \cite{Kato78}.  In the following we will prove that a compact component in
the moduli space would be preserved under small deformations.  We do this using a third guise of our moduli space, justifying at
the same time, finally, why we speak of ``$\PU(2)$-instantons''.

Let $P$ be the principal $\PU(2)$-bundle obtained as the quotient of the principal $\U(2)$ frame bundle of $E$ by the centre of
$\U(2)$.  Remark that the adjoint action $\Ad$ of $\SU(2)$ on itself descends to an action of $\PU(2)\isom\SU(2)/\{\pm1\}$ on
$\SU(2)$ so that we can define the gauge group $\mathscr G\coloneq\Gamma(P\times_{\Ad}\SU(2))$.  This group acts naturally on
the affine space $\mathscr A$ of connections on $P$.  We call a connection irreducible if its stabiliser in $\mathscr G$ is
minimal, i.\,e.\ the center $\{\pm1\}$ of $\mathscr G$, and denote by $\mathscr A^*$ the space of irreducible connections.  The
moduli space of irreducible anti-self dual connections on $P$ is now defined as the quotient
\[
	\mathcal M^{\rm ASD}(P)^*\coloneq\{A\in\mathscr A^*\colon F_A^+=0\}/\mathscr G
\]
where $F_A^+$ denotes the self-dual part of the curvature $F_A$ of $A$.  There is a canonical isomorphism
\begin{equation}
	\label{eq:PE}
	\mathcal M^{\rm ASD}(E)^*\isom\mathcal M^{\rm ASD}(P)^*
\end{equation}
with the moduli space of irreducible anti-self-dual connections on $E$ from the previous section, independent of the fixed
connection $a$ on $\det E$.  This independence will allow us to construct a parametrised moduli space for a deformation of our
surface.

To do this we write this moduli space in a different way as follows.  The space $\mathscr A^*$ is a principal $\mathscr
G/\{\pm1\}$-bundle over the corresponding orbit space $\mathscr B^*\coloneq\mathscr A^*/\mathscr G$.  The map $F^+\colon\mathscr
A\to\Omega^2_+(\ad P)$ associating to a connection $A$ the self-dual part $F_A^+$ of its curvature is $\mathscr G$-equivariant
and therefore defines a section $F^+\colon\mathscr B^*\to\mathscr E$ in the associated vector bundle $\mathscr E\coloneq\mathscr
A^*\times_{\ad}\Omega^2_+(\ad P)$ over $\mathscr B^*$.  The moduli space $\mathcal M^{\rm ASD}(P)^*$ is then simply the
vanishing locus of this section.  Using suitable \name{Sobolev} completions $F^+$ is a \name{Fredholm} map between \name{Banach}
manifolds.  A set $\mathcal C\subset\mathcal M^{\rm ASD}(P)^*$ is said to be regular if $F^+$ is regular at every point of
$\mathcal C$.  The following proposition allows one to check regularity using the complex geometric framework.  It results from
comparing the local models of the moduli spaces $(\mathcal M^{\rm ASD})^*$ and $\mathcal M^{\rm st}$.

\begin{proposition}[\cite{Luebke&Teleman}]
	A point in $(\mathcal M^{\rm ASD})^*$ is regular if and only if its image in $\mathcal M^{\rm st}$ under the
	\name{Kobayashi}-\name{Hitchin}-correspondence \eqref{eq:KH} is regular.
\end{proposition}

\begin{corollary}
	\label{cor:regular}
	For a class $\VII_0^1$ surface $S$ every compact component $\mathcal C\subset\mathcal M^{\rm ASD}(S)^*$ is regular.
\end{corollary}

\begin{proof}
	By proposition~\ref{prop:local}, $\mathcal M^{\rm st}(S)$ is regular at every smooth point and we saw that all singular
	points lie on a non-compact component.
\end{proof}

We show that in general a regular compact component of the moduli space of irreducible ASD connections is stable under small
deformations of the metric.  For this we consider a parametrised version of the above construction of the moduli space $\mathcal
M^{\rm ASD}(P)^*$.  Let $I$ be the interval $[-1,+1]$ and $(g_t)_{t\in I}$ a smooth one-parameter family of \name{Riemann}ian
metrics $g_t$ on the base manifold.  Again, $\underline{\mathscr A}^*\coloneq\mathscr A^*\times I$ is a principal $\mathscr
G/\{\pm1\}$-bundle over $\underline{\mathscr B}^*\coloneq\mathscr B^*\times I$.  The map $\underline
F^+\colon\underline{\mathscr A}\to\Omega^2(\ad P)$, defined by mapping $(A,t)$ to the self-dual part $F_A^{+_{g_t}}$ of the
curvature $F_A$ with respect to the metric $g_t$, is $\mathscr G$-equivariant and defines a section $\underline{\mathscr
B}^*\to\underline{\mathscr A}^*\times_{\ad}\Omega^2(\ad P)$.  This section actually takes values in the subbundle
$\underline{\mathscr E}$ whose fibre over $([A],t)$ is the space $\Omega^2_{+_{g_t}}(\ad P)$ of $(\ad P)$-valued two-forms that
are self-dual with respect to the metric $g_t$.  This gives a section $\underline F^+\colon\underline{\mathscr
B}^*\to\underline{\mathscr E}$ whose vanishing locus is the parametrised moduli space
\[
	\big(\underline{\mathcal M}^{\rm ASD}\big)^*\coloneq\{([A],t)\in\mathscr B^*\times I\colon F_A^{+_{g_t}}=0\}
	\fullstop
\]
The restriction $\pi\colon\big(\underline{\mathcal M}^{\rm ASD}\big)^*\to I$ of the projection $\mathscr B^*\times I\to I$
gives a fibration
\[
	\big(\underline{\mathcal M}^{\rm ASD}\big)^*=\bigcup_{t\in I}\pi^{-1}(t)
	\qquad
	\text{with}
	\quad
	\pi^{-1}(t)=\mathcal M^{\rm ASD}(g_t)^*\times\{t\}
	\fullstop
\]

\begin{proposition}
	For $t$ sufficiently small $\mathcal M^{\rm ASD}(g_t)^*$ contains a regular compact component if $\mathcal
	M^{\rm ASD}(g_0)^*$ does.
\end{proposition}

\begin{proof}
	Let $\mathcal C\subset\mathcal M^{\rm ASD}(g_0)^*$ be such a regular compact component.  The restriction of $\underline
	F^+$ to $\mathcal M^{\rm ASD}(g_0)^*=\pi^{-1}(0)$ is just the above map $F^+$ and thus regular on $\mathcal C$.  Therefore
	$\underline F^+$ itself is regular on $\mathcal C$.  Regularity is an open condition so $\underline F^+$ is regular on an
	open neighbourhood $N$ of $\mathcal C$ in $\big(\underline{\mathcal M}^{\rm ASD}\big)^*$.  It follows that $N$ is a
	finite-dimensional smooth open manifold.  Then, as $\mathcal C$ is compact, we can choose a compact neighbourhood $K$ of
	$\mathcal C$ in $N$ with $K\cap\pi^{-1}(0)=\mathcal C\subset\mathring K$.  We have $\mathring
	K\cap\pi^{-1}(0)=K\cap\pi^{-1}(0)$.  It suffices to show that $\mathring K\cap\pi^{-1}(t)=K\cap\pi^{-1}(t)$ for $t$
	sufficiently small.  Suppose not.  Then there exists a sequence of points $([A_n],t_n)\in(K\setminus\mathring
	K)\cap\pi^{-1}(t_n)$ with $t_n\to0$.  But $K$ being compact, some subsequence of it converges to a point
	$([A],0)\in(K\setminus\mathring K)\cap\pi^{-1}(0)=\emptyset$ which is a contradiction.
\end{proof}

\begin{corollary}
	\label{cor:connected}
	For a class $\VII_0^1$ surface $S$, all moduli spaces $\mathcal M^{\rm ASD}(S)^*\isom\mathcal M^{\rm st}(S)$, $\mathcal
	M^{\rm pst}(S)$ and $\mathcal M^{\rm s}(S)$ are connected.
	\hfill\qedsymbol
\end{corollary}

\begin{proof}
	We saw that another connected component in one of these moduli spaces, other than the one we already described, would belong
	to $\mathcal M^{\rm ASD}(S)^*$ and therefore be compact.  By corollary~\ref{cor:regular} it
	would also be regular.  Let now $(J_t)_{t\in I}$ be a family of complex structures on the real manifold underlying $S$,
	parametrising a degeneration $(S_t)_{t\in I}$ of blown-up primary \name{Hopf} surfaces $S_t$, $t\not=0$, into $S_0\coloneq
	S$.  We can take $(g_t)_{t\in I}$ to be a corresponding smooth family of \name{Gauduchon} metrics $g_t$ on $S_t$.  Then the
	preceding corollary says that $\mathcal M^{\rm ASD}(S_t)^*$ would contain a compact component too, contradicting
	\cite{Toma06}.
\end{proof}

\section{The moduli spaces}
\label{sec:results}

We can finally assemble all our results to a complete description of the moduli spaces.  By a \dfn{compact complex space with
smooth boundary} we mean a compact real analytic space with a smooth boundary structure around its boundary and a possibly
singular complex structure on its interior.  We write ``$\lle$'' for ``$<(\le)$''.

\begin{theorem}
	Let $S$ be a minimal class \VII surface with $b_2(S)=1$.
	\newcounter{contenumi}
	\begin{enumerate}
		\item \label{enum:<0}
			If $\deg K<0$ --- i.\,e.\ if $S$ is the half or the parabolic \name{Inoue} surface or an \name{Enoki} surface with
			$\deg K<0$ --- then the entire moduli space $\mathcal M^{\rm (p)st}(S)$ of (poly)stable holomorphic bundles
			of type \eqref{eq:fix} is bijectively parametrised by the open (closed) complex one-dimensional disc
			$\Pic^0_{\lle\varrho}(S)\cup\{0\}$.
		\item \label{enum:>0}
			If $\deg K\ge0$ --- i.\,e.\ $S$ is an \name{Enoki} surface with $\deg K\ge0$ --- then the $\mathcal
			M^{\rm (p)st}(S)$ is bijectively parametrised by
			$\big(\!\Pic^0_{\lle\varrho}(S)\cup\{0\}\big)\setminus\:\unstable(S)$ where $\unstable(S)$ is the finite set
			\eqref{eq:unstables}.
		\setcounter{contenumi}{\value{enumi}}
	\end{enumerate}
	The parametrisation is given by mapping
	\[
		\Pic^0_{=\varrho}(S)\ni\hol L\mapsto\hol L\oplus(\hol L^\dual\otimes\hol K)
		\comma
		\qquad
		\Pic^0_{<\varrho}(S)\ni\hol L\mapsto\hol E_{\hol L}
		\qquad\text{and}\qquad
		0\mapsto\hol E
		\comma
	\]
	where:
	\begin{enumerate}
		\setcounter{enumi}{\value{contenumi}}
		\item \label{enum:filtrable}
			On the half \name{Inoue} surface, $\hol E$ is the filtrable bundle $\hol A_{\hol O}$ and $\mathcal M^{\rm (p)st}(S)$
			contains no non-filtrable bundles.
		\item \label{enum:non-filtrable}
			On an \name{Enoki} or parabolic \name{Inoue} surface $\hol E$ is the only non-filtrable bundle in $\mathcal
			M^{\rm (p)st}(S)$.
	\end{enumerate}
	In case (\ref{enum:<0}) this is a homeomorphism, holomorphic on the stable part.
	In case (\ref{enum:>0}) this is a local homeomorphism except at points $\hol R\in\roots(S)$, holomorphic on the stable part
	minus $\roots(S)$.  $\mathcal M^{\rm st}(S)$ is a one-dimensional complex space whose singularities are simple normal
	crossings at the points $\hol E_{\hol R}$ characterised by
	\[
		 \lim_{\hol L\to\hol R^\dual(-C)}\hol E_{\hol L}
		=\hol E_{\hol R}
		=\lim_{\hol L\to\hol R          }\hol E_{\hol L}
		\qquad
		\text{for}
		\quad
		\hol R^\dual(-C)\in\unstable(S)
		\fullstop
	\]
	Their number $|\unstable(S)|$ is finite but unbounded if the metric varies in the space of \name{Gauduchon} metrics.

	Therefore, except for the case $\Pic^0_{=\varrho}(S)\cap\roots(S)\not=\emptyset$ on an \name{Enoki} surface, $\mathcal
	M^{\rm pst}(S)$ is a one-dimensional compact complex space with smooth boundary a circle and interior $\mathcal
	M^{\rm st}(S)$, smooth in case (\ref{enum:<0}) and in general singular in case (\ref{enum:>0}).
\end{theorem}

\begin{figure}[h]
	\centering
	\includegraphics[width=\textwidth]{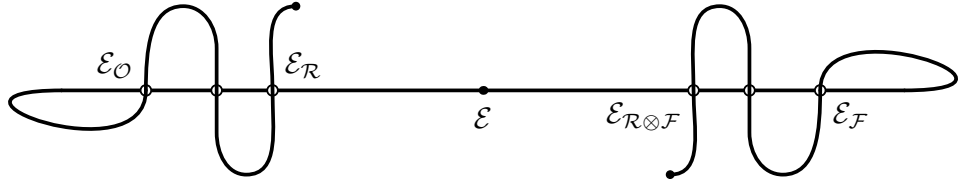}
	\caption{The moduli space of polystable bundles for an \name{Enoki} surface}
	\label{fig:poly}
\end{figure}

For an \name{Enoki} surface $S$ the moduli space $\mathcal M^{\rm pst}(S)$ can be viewed as a closed complex disc with finitely
many self intersections as in figure~\ref{fig:poly} (where we divided dimensions by two).  Note that the degree corresponds to
(the logarithm of) the ``distance'' from the center of the disc.  In the limit case where a line bundle $\hol R\in\roots(S)$
happens to lie on the boundary circle of this disc, the self intersection is merely a ``touch'' of a point on the boundary
circle with an interior point, but both points do not belong to the moduli space since they correspond to the unstable bundles
$\hol E_{\hol R}$ and $\hol E_{\hol R^\dual(-C)}$.

Since non-filtrable bundles are stable by definition, the above also completes our description of the moduli space $\mathcal
M^{\rm s}(S)$ of simple holomorphic bundles of type \eqref{eq:fix}.  If $S$ is the parabolic \name{Inoue} surface then
$\mathcal M^{\rm s}(S)$ is simply isomorphic to $(\Pic^0(S)\cup\{0\})\setminus\queer(S)$, i.\,e.\ to the complex line $\C$
minus a discrete set of points.

\begin{figure}[h]
	\centering
	\includegraphics{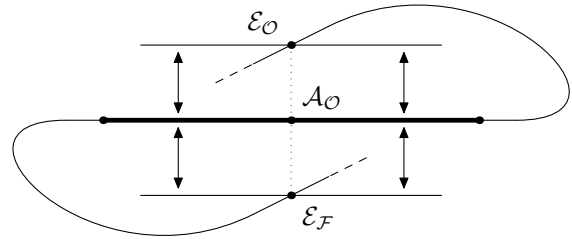}
	\caption{The moduli space for the half \name{Inoue} surface}
	\label{fig:half}
\end{figure}

If $S$ is the half \name{Inoue} surface then, due to the isomorphism $\hol A_{\hol F}\isom\hol A_{\hol O}$, the smooth branches
in the two local pictures in figure~\ref{fig:local} for $\hol R=\hol O$ and $\hol R=\hol F$ coincide.  With notations as in
proposition~\ref{prop:local}, we can regard the curves $C''_{\hol O}=C''_{\hol F}$, $C'_{\hol O}$ and $C'_{\hol F}$ as one
single ``curve'' with a triple point consisting of the three non-separable points $\hol A_{\hol O}$, $\hol E_{\hol O}$ and $\hol
E_{\hol F}$.  This curve is smooth at $\hol A_{\hol O}$ but transversely crossed by $C_{\hol O}$ at $\hol E_{\hol O}$ and by
$C_{\hol F}$ at $\hol E_{\hol F}$.  The resulting moduli space $\mathcal M^{\rm s}(S)$ is depicted in figure~\ref{fig:half}
(where the stable part is marked in bold and we omitted indicating the punctures corresponding to bundles in $\queer(S)$).

\begin{figure}[h]
	\centering
	\includegraphics[width=\textwidth]{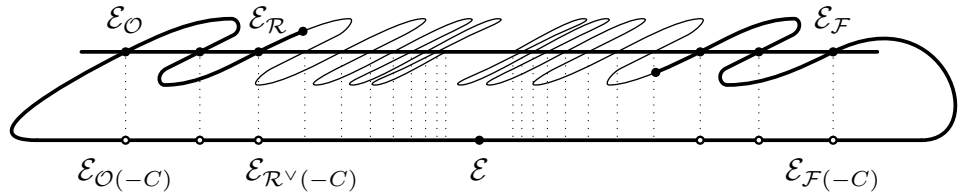}
	\caption{The moduli space for an \name{Enoki} surface}
	\label{fig:simple}
\end{figure}

If $S$ is an \name{Enoki} surface then $\mathcal M^{\rm s}(S)$ contains no such triple points but countably infinitely many
pairs of inseparable points $\hol E_{\hol R}$ and $\hol E_{\hol R^\dual(-C)}$ corresponding to line bundles $\hol
R\in\roots(S)$, the first of them being singular and the second smooth as in figure~\ref{fig:local}.  This is shown in
figure~\ref{fig:simple}.

\bibliographystyle{amsalpha}
\bibliography{\jobname}

\end{document}